\documentclass[11pt,centertags]{article}

\usepackage{amsmath}
\usepackage{amsthm}
\usepackage{amscd}
\usepackage{amssymb}
\usepackage{verbatim}
\usepackage[all]{xy}
\usepackage{latexsym}
\usepackage[dvips]{graphics}
\usepackage{psfrag}

\title{Some recent applications of the barycenter method in geometry}
\author{Christopher Connell and Benson Farb~\thanks{Both authors are
supported in part by the NSF.  This paper grew out of talks given by the
first and second authors at the Cornell and Georgia Topology
conferences, respectively.}\\
Department of Mathematics\\
University of Chicago
}
\theoremstyle{definition}
\newtheorem*{Definition}{Definition}

\theoremstyle{plain}
\newtheorem{theorem}{Theorem}
\newtheorem{proposition}[theorem]{Proposition}
\newtheorem{lemma}[theorem]{Lemma}
\newtheorem{corollary}[theorem]{Corollary}

\newtheorem{conjecture}[theorem]{Conjecture}
\newtheorem{question}[theorem]{Question}
\newtheorem{prob}[theorem]{Problem}
\newtheorem{xample}[theorem]{Example}

\newcounter{remarks}
\newenvironment{remarks}
{\paragraph*{Remarks}\smallskip
     \begin{list}{\arabic{remarks}. }{\usecounter{remarks}
          \setlength{\leftmargin}{0in}
          \setlength{\rightmargin}{0in}
          \setlength{\labelsep}{0pt}
          \setlength{\labelwidth}{0pt}
          \setlength{\listparindent}{0pt}
     }
}
{
\end{list}
}

\def\eps{\epsilon}

\def\bar{\overline}

\newcommand\D{\partial}

\newcommand\op{\operatorname}
\newcommand\ds{\displaystyle}
\newcommand\Jac{\operatorname{Jac}}

\renewcommand\deg{{\rm deg}}
\renewcommand\tilde{\widetilde}
\newcommand\vol{\operatorname{Vol}}
\newcommand\minvol{\operatorname{Minvol}}
\newcommand\Minvol{\operatorname{Minvol}}
\newcommand\directsum{\oplus}
\newcommand\htop{\op{h}_{\op{top}}}
\newcommand\R{\mbox{\bf R}}
\newcommand\hyp{\mbox{\bf H}}

\newcommand\Rrank{\R\mbox{-rank}}

\DeclareMathOperator{\Vol}{Vol}
\DeclareMathOperator{\Isom}{Isom}
\DeclareMathOperator{\ent}{ent}
\DeclareMathOperator{\Id}{Id}
\DeclareMathOperator{\meas}{{\cal M}}
\DeclareMathOperator{\bary}{bar}
\DeclareMathOperator{\supp}{supp}
\DeclareMathOperator{\SO}{SO}
\DeclareMathOperator{\rank}{rank}
\DeclareMathOperator{\SL}{SL}
\renewcommand{\to}{\longrightarrow}
\newcommand{\Mvariable}[1]{\, #1 \,}

\begin{document}
\maketitle
\tableofcontents

\pagebreak

\section{The Entropy Rigidity Conjecture}
In this paper we describe some recent applications of the barycenter
method in geometry.  This method was first used by Duady-Earle and later
greatly extended by Besson-Courtois-Gallot in their solution of a number
of long-standing problems, in particular in their proof of
entropy rigidity for closed, negatively curved locally symmetric
manifolds.  Since there are already a number of surveys describing
this work (see \cite{BCG2,BCG3,Pa,Ga}), we will concentrate here only on
advances that have occured after these surveys appeared.  While most of
this paper is a report on results appearing in other papers,
some of the material here is new (e.g. Proposition
~\ref{proposition:convexity} and the examples in Section 
~\ref{section:counterexamples}).

\subsection{The main conjecture}
The {\em volume entropy} of a closed Riemannian $n$-manifold $(M,g)$,
denoted by $h(g)$, is defined to be

$$h(g)=\lim_{R\rightarrow \infty}\frac{1}{R} \log (\Vol(B(x,R)))$$
where $B(x,R)$ is the ball of radius $R$ around a fixed point $x$ in
the universal cover $X$.  The number $h(g)$ is independent of the
choice of $x$, and equals the topological entropy of the geodesic flow
on $(M,g)$ when the curvature $K(g)$ satisfies $K(g)\leq 0$ (see
\cite{Ma}).  Note that while neither the volume $\Vol(M,g)$ nor the
entropy $h(g)$ is invariant
under scaling the metric $g$, the {\em normalized entropy}
$$\ent(g)=h(g)^n\Vol(M,g)$$ is scale
invariant.

Now let $M$ be a close $n$-manifold that admits a locally symmetric
Riemannian metric $g_{loc}$ of nonpositive sectional curvature.  When
$(M,g_{loc})$ is not locally isometric to a product then $g_{loc}$ is
unique up to homothety (i.e. multiplying the metric by a number).  When
$M$ is locally a product, one may show using Lagrange multipliers that
there is a unique (up to homothety) locally symmetric metric on $M$
which minimizes $\ent(g)$; see \cite{CF1}.  Henceforth we will abuse
notation and denote this
metric by $g_{loc}$, and call it {\em the} locally symmetric
metric.

The {\em Entropy Rigidity Conjecture}, stated in various forms
by Katok, Gromov, and Besson-Courtois-Gallot (see, e.g., \cite{BCG2},
Open Question 5), has two components.  The purely metric component
posits that for most $M$, the metric $g_{loc}$ minimizes the functional
$\ent(g)$ over the space of all Riemannian metrics on $M$, and in fact
{\em uniquely} minimizes $\ent(g)$ (up to homothety).  This property
would characterize the locally symmetric metric by essentially a
single number.  The topological component of the conjecture is an
extension of this statement from manifolds to maps.

\begin{conjecture}[Entropy Rigidity Conjecture]
\label{conjecture:entrig}
Let $M$ be a closed manifold which admits a locally symmetric Riemannian
metric $g_{loc}$ with nonpositive sectional 
curvature. Assume that $(M,g_{loc})$ has
no local factors isometric to $\R$.  Let $(N,g)$ be any
closed Riemannian manifold, and let $f:N\to M$ be any continuous map.
Then
\begin{equation}
\label{eq:entrig}
\ent(N,g)\geq |\deg f|\ent(M,g_{loc})
\end{equation}
with equality if and only if $f$ is
homothetic to a Riemannian covering.
\end{conjecture}

\bigskip
\noindent
{\bf Remarks: }
\begin{enumerate}
\item  The case when $f$ is a homeomorphism, or even the identity map,
gives that $g_{loc}$ uniquely minimizes $\ent(g)$.
\item Easy examples show that the
restrictions to $M$ with no local $\R$ factors is necessary.
\item Conjecture~\ref{conjecture:entrig} easily implies Mostow Rigidity;
the argument is the same as in the case when $M$ is negatively curved,
which is given in \cite{BCG2}.
\end{enumerate}
Conjecture~\ref{conjecture:entrig} may be extended to a variety
of other contexts, for example to finite volume metrics, Finsler
metrics, foliations with locally symmetric leaves, contact flows,
and magnetic flows. To be more precise, the table below shows
some natural choices for the volume and entropy in various metric
settings. In the table we have adopted the following conventions:

\begin{itemize}
\item $\theta$ represents a contact $1$-form
\item $h_m$ is metric entropy, i.e. the measure 
theoretic entopy for the Liouville
measure
\item $h_{\op{top}}$ represents the topological entropy of the geodesic flow
\item $h_{\op{vol}}$  is the volume growth entropy
\item $\dim_{\mathcal{H}}(\Lambda(\Gamma))$  is the Hausdorff dimenon of
the limit set $\Lambda(\Gamma)$ of a discrete group $\Gamma$ of
hyperbolic isometries
\item  $M=Y/\Gamma$
\end{itemize}

\begin{figure}[ht]
\label{table:settings}

\begin{center}

\setlength{\baselineskip}{8pt} \small
\begin{tabular}{|p{3cm}|p{3cm}|p{5.5cm}|} \hline
{\bf Setting}  & {\bf Volume $\op{Vol}(M)$} &  {\bf entropy $h(d)$}\\
\hline\hline
\parbox[3cm]{6cm}
{Riemannian metrics} on closed $M$& Riem. volume & $h_{\op{vol}}$  \\
\hline Finite vol. metrics on noncompact $M$ & Riem. volume 
& $\inf\{s : \int_Y e^{-s d(p,y)}dg(y)
<\infty\}$  \\ 

\hline Contact flows &
$\ds\int_M\theta\wedge(d\theta)^{n-1}$ & $h_{\op{top}}$ \\ \hline
Finsler metrics & Finsler Volume & $\htop$\\
\hline Magnetic field flows & Euler-Lagrange Volume & $h_{\op{top}}$ \\
\hline Foliations on $M$ &
$d\op{Vol}_{g_L}\times d\op{Vol}(L)$: leafwise volume $\times$ transverse
volume &
$\ds\left(\frac{1}{\op{Vol}(M)}\int_M h(g_L)^nd\op{Vol}\right)^{1/n}$ \\
\hline
Geometrically finite $M$ & 
Volume of core & $\dim_{\mathcal{H}}(\Lambda(\Gamma))$  \\
\hline
\end{tabular}

\end{center}

\caption{Some natural choices for the volume and entropy in different
settings.  The functional being minimized is 
$\op{ent}:=h(d)^n\op{Vol}(M)$.}
\end{figure}

These extensions are discussed in greater depth in Sections
~\ref{section:negative} and~\ref{section:positive}.

\subsection{Progress}

The main evidence for Conjecture~\ref{conjecture:entrig} is given
by the following theorem.

\begin{theorem}[Besson-Courtois-Gallot \cite{BCG1}]
\label{theorem:bcg}
Conjecture~\ref{conjecture:entrig} is true when $(M,g_{loc})$ has
$\Rrank$ one, i.e. when $(M,g_{loc})$ is negatively curved.
\end{theorem}

This result and its proof have a number of
corollaries, including solutions to
long-standing problems on geodesic flows, asymptotic harmonicity, and
Gromov's minvol invariant; these are described in \cite{BCG2}.

In higher rank very little is known.  The following
was announced in \cite{BCG2} and later in \cite{BCG3}, and was proved
in \cite{CF1} (in the finite volume case as well).

\begin{theorem}
Conjecture~\ref{conjecture:entrig} is true when $(M,g_{loc})$ is
locally (but not necessarily globally)
isometric to a product of negatively curved locally symmetric spaces.
\end{theorem}

It seems that a significantly new idea is needed to prove Conjecture
~\ref{conjecture:entrig}.  See, however, \S~\ref{section:degree} below.

\subsection{The Besson-Courtois-Gallot map}
\label{section:degree:argument}

Building on earlier work of Duady-Earle, Besson-Courtois-Gallot
constructed a remarkable map $F$ in every homotopy class of maps
$f:N\rightarrow M$ from an arbitrary closed manifold $N$ to any closed
manifold $M$ of negative sectional curvature.  This map is a kind of
replacement for, and has several advantages over, the harmonic map in
the homotopy class of $f$; in particular, it can often be differentiated
explicitly.  The map $F$ is a key ingredient in the proof of Theorem
~\ref{theorem:bcg}.

Now suppose that, in addition to being negatively curved, 
$M$ is also locally symmetric with metric $g_{loc}$.
In this case Besson-Courtois-Gallot were able to give a precise estimate
on the Jacobian $|\Jac F|$, proving Theorem~\ref{theorem:bcg}.  To state
things precisely, it will be useful to fix a parameter $s>h(g)$, where $g$
is the metric on $N$ and $h(g)$ is the volume growth entropy of
$g$.  When $(M,g_{loc})$ has $\Rrank$ one, that is
when it is negatively curved, Besson-Courtois-Gallot proved
that, given any continuous $f:N\to M$, there
exists a smooth map $F_s:N\to M$ homotopic to $f$ with the following two
important properties:

\bigskip
\noindent
{\bf Universal Jacobian bound: } For all $s>h(g)$ and all $y\in N$ we have:
\begin{equation}
\label{equation:jacobian:bound}
\vert \Jac F_s(y)\vert\le \left(\frac s{h(g_{loc})}\right)^n
\end{equation}

\bigskip
\noindent
{\bf Infinitessimal rigidity: } 
\label{equation:infinitessimal}
There is equality in equation
(~\ref{equation:jacobian:bound}) at the point $y\in N$ if and only if
$D_yF_s$ is a homothety.
\bigskip

Before explaining how to construct such an $F_s$ (which we do in
\S~\ref{section:construction}), let us deduce some consequences of its
existence.

The existence of any $F$ satisfying the two properties above implies
Conjecture~\ref{conjecture:entrig} by an elementary degree argument, as
follows.  Since for $s>h(g)$, the map $F_s$ is a $C^1$ map, we may
simply compute:

\begin{eqnarray}
|\deg(f)| \vol(M)&=\left| \int_{N}f^*dg_{loc}
\right| \\
& \\
&  =\left| \int_{N} F_s^*dg_{loc} \right|\\
& \\
& \leq\int_{N} \left| \Jac F_s\right| dg \\
& \\
&\leq
\left(\frac{s}{h(g_{loc})}\right)^n \vol(N)
\end{eqnarray}

Letting $s\to h(g)$ gives the inequality in Conjecture
~\ref{conjecture:entrig}. In the case when equality is achieved,
after scaling the metric $g$ by the constant
$\frac{h(g)}{h(g_{loc})}$, we have $h(g)=h(g_{loc})$ and
$\Vol(N)=|\deg(f)|\,\Vol(M)$.

\subsection{The degree theorem}
\label{section:degree}

It is useful for us to restate the inequality (~\ref{eq:entrig}) of the
Entropy Rigidity Conjecture as
\begin{equation}
\label{equation:degree}
\deg(f)\leq C\frac{\vol(N)}{\vol(M)}
\end{equation}
where $C=\left(\frac{h(g)}{h(g_{0})}\right)^n$.  In particular, proving the
inequality (~\ref{equation:degree}) with this precise value of $C$
gives that the locally symmetric metric $g_{loc}$ minimizes $\ent(g)$
among all metrics $g$ on $N$.

When the inequality (~\ref{equation:degree}) holds for {\em some}
(universal) $C$, perhaps bigger than
$\left(\frac{h(g)}{h(g_{0})}\right)^n$, then the degree argument given
in \S~\ref{section:degree:argument} gives a universal bound
on the degree of any map of any manifold into $N$.  This can be viewed
as the topological part of Conjecture~\ref{conjecture:entrig}.  When
$N$ is negatively curved, Gromov
proved this fundamental relationship between degree and volume in his
paper \cite{Gr} (see his Volume Comparison Theorem).

We note that by {\em universal} we mean that $C$ depends only on the
dimension of the manifolds $M,N$ and on their smallest Ricci curvatures.
Scaling the metrics shows that these dependencies are necessary.

When the sectional curvatures of $N$ are not necessarily negative
but only nonpositive, the situation is more complicated.
Since the $n$-dimensional torus has flat metrics and also has self-maps
of arbitrary degree, one needs at least to assume that $N$ has no local
$\R$ factors in order that (~\ref{equation:degree}) be true.  
Using the barycenter technique (see
below), we proved the following in \cite{CF2}.

\begin{theorem}[The Degree Theorem]
\label{theorem:degree}
Let $M$ be a closed, locally symmetric $n$-manifold with nonpositive
sectional curvatures.  Assume that $M$ has no local direct factors
locally isometric to $\R, \hyp^2$, or $\SL_3(\R)/\SO_3(\R)$.  Then for
any closed Riemannian manifold $N$ and any continuous map $f:N\to M$,
$$\deg(f)\leq C\frac{\vol(N)}{\vol(M)}$$ where $C>0$ depends only on $n$
and the smallest Ricci curvatures of $N$ and $M$.
\end{theorem}

While the ``no local $\R$ factors'' assumption is necessary, we believe
that Theorem~\ref{theorem:degree} should hold in all other cases.  The
case when $M$ is locally modelled on $\SL_n(\R)/SO_n(\R), n\geq 2$ also
follows from work of Savage \cite{Sa}; see \S~\ref{section:gromov:norm} below.

Theorem~\ref{theorem:degree} follows immediately
from the degree argument given above (in
\S~\ref{section:degree:argument}) together with the following theorem,
which is the main theorem of \cite{CF2}, and which
we believe is of independent interest.

\begin{theorem}[Universal Jacobian bound]
Any continuous map $f:N\rightarrow M$ between closed
$n$-manifolds, with $M$ nonpositively curved and locally symmetric (barring
the exceptions of Theorem~\ref{theorem:degree}),
is homotopic to a $C^1$ map $F$ with
$$|\Jac F|\leq C$$
for some constant $C$
depending only on $n$ and on the smallest Ricci curvatures of
$M$ and $N$.
\end{theorem}

The map $F$ is the extension to nonpositively curved locally symmetric
manifolds of the Besson-Courtois-Gallot map; we discuss this below.

As we will discuss in \S~\ref{section:jacobian:bound},
the universal bound for $|\Jac F|$
is obtained by reducing the problem to a minimization problem over a
space of measures on a certain Lie group.  The Lie group appears because the
natural boundary attached to symmetric spaces may be described {\em
algebraically}.  For general nonpositively curved manifolds no such
description is available.  But the following still seems possible.

\begin{conjecture}[Degree Conjecture in nonpositive curvature]
\label{conjecture:degree}
Let $M$ be a closed $n$-manifold with nonpositive sectional curvature,
negative Ricci curvature, and no local $\R$ factors.  Then for any
closed Riemannian manifold $N$ and any continuous map $f:N\to M$,
$$\deg(f)\leq C\frac{\vol(N)}{\vol(M)}$$ where $C>0$ depends only on $n$
and the smallest Ricci curvatures of $N$ and $M$.
\end{conjecture}

The evidence for Conjecture~\ref{conjecture:degree} is that it is
true when $M$ has negative sectional curvature (Gromov \cite{Gr}), and
(almost always) when $M$ is locally symmetric
(by Theorem~\ref{theorem:degree}, since
no local $\R$ factors implies negative Ricci curvature for locally
symmetric manifolds of nonpositive sectional curvature).  The
``negative Ricci curvature'' hypothesis in Conjecture
~\ref{conjecture:degree} is necessary, as the following example
(see \cite{BGS}) shows.

\begin{xample}
\label{example:gluing}
{\rm 
Let $S$ be a negatively curved surface with one boundary component
$C$. While keeping $S$ nonpositively curved we may smooth a neighborhood
of $C$ into a cylinder which is metrically $[0,\eps)\times S^1$, where
$C$ is identified to $\{0\}\times S^1$. Call the resulting surface
$S'$. 

Now take two copies $X_1$ and $X_2$ of the manifold $S'\times S^1$
and glue $X_1$ to $X_2$ by identifying $[0,\eps/3]\times S^1\times
S^1\subset X_1$ and $[0,\eps/3]\times S^1\times S^1\subset X_2$ by the
map $(t,x,y)\mapsto(\eps/3-t,y,x)$.  The resulting $3$-manifold $X$
has no flats, but has some zero curvature at every point and even has an
open subset of zero curvature, namely the subset $Z\subset X$ formed
from the union of the images of $[0,\eps)\times S^1$ in each $X_i$. 

We note that the two overlapping pieces, $(X_1\setminus Z)\cup
(\eps/2,\eps)\times S^1\times S^1$ and $(X_1\setminus Z)\cup
(\eps/2,\eps)\times S^1\times S^1$, each carry an $S^1$ action that is
coherent with a torus action on $Z$ which acts the same in each
factor. Therefore, the whole of $X$ carries an {\em $F$-structure} in
the sense of Cheeger and Gromov \cite{CheegerGromov}.  This implies via
\cite{Paternain,CheegerGromov} that $\op{Minvol}(X)=0$ (see below for
the definition of $\op{Minvol}(X)$); in
particular, the conclusion of Conjecture~\ref{conjecture:degree} 
cannot hold for $X$, even when $f$ is the identity map.  

Note that $X$ has no local $\R$ direct factors, but that $X$ does have
some points with some direction of zero Ricci curvature.}
\end{xample}

\subsection{A related conjecture}
\label{section:gromov:norm}

Gromov has made the following:

\begin{conjecture}[Positivity of Gromov norm]
\label{conjecture:gromov:norm}
The Gromov norm of a closed, nonpositively curved, locally symmetric
manifold $M$ with no local $\R$ factors is positive.
\end{conjecture}

Conjecture~\ref{conjecture:gromov:norm} was proven by Savage \cite{Sa}
when $M$ is locally isometric to $\SL_n(\R)/SO_n(\R)$. 
While we do not see how positivity of Gromov
norm for $M$ directly implies the inequality
(~\ref{equation:degree}), Savage's proof does imply this
inequality.  The key point is that Savage proves that any simplex
in $M$ can be ``straightened'' to a simplex with universally
bounded volume.  This should hold for all symmetric spaces of
noncompact type, but this is still an open question.

\subsection{Some consequences of the Degree Theorem}

We end this section by recalling some consequences, given in \cite{CF2},
of Theorem~\ref{theorem:degree}.

\bigskip
\noindent
{\bf The Minvol invariant.} One of the basic invariants associated to a
smooth manifold $M$ is its {\em minimal volume}:
$$\minvol(M):=\inf_g\{\Vol(M,g):|K(g)|\leq 1\}$$ where $g$ ranges over
all smooth metrics on $M$ and $K(g)$ denotes the sectional curvature of
$g$.  The basic questions about $\minvol(M)$ are: for which $M$ is
$\minvol(M)>0$? when is $\minvol(M)$ realized by some metric $g$?

When a nonpositively curved manifold $M$ has a local direct factor
locally isometric to $\R$, then $M$ has some finite cover $M'$ with an $S^1$
direct factor, which implies that $M'$ has positive degree self-maps and
that $\minvol(M')=0$.  On the other hand, by taking $f$ in Theorem
~\ref{theorem:degree} to be the identity map while allowing the metric
$g$ on $M$ to vary, immediately gives that $\minvol(M)>0$ for locally
symmetric $M$, barring the exceptional cases in Theorem
~\ref{theorem:degree}.  However, positivity was already known more
generally, by the following result of Gromov (see also \cite{Sa} for the
case of $M$ locally isometric to the symmetric space for $\SL(n,\R)$).

\begin{corollary}[Positivity of Minvol]
\label{corollary:minvol}
Let $M$ be a closed, locally symmetric manifold with nonpositive
curvature and no local $\R$ factors.  Then $\minvol(M)>0.$
\end{corollary}

When $M$ is real hyperbolic, Besson-Courtois-Gallot \cite{BCG1} proved
that $\minvol(M)$ is {\em uniquely} realized by the locally symmetric
metric.

\begin{prob}
Compute $\minvol(M)$ for all closed $M$ which admit a locally symmetric
metric of nonpositive curvature, with no local $\R$ factors.
Is $\minvol(M)$ always realized by the locally symmetric metric? Is it
realized uniquely?
\end{prob}

\bigskip
\noindent {\bf Self maps and the co-Hopf property. }
As $\deg(f^n)=\deg(f)^n$, an immediate
corollary of Theorem~\ref{theorem:degree} is the following.

\begin{corollary}[Self maps]
\label{corollary:hopf}
Let $M$ be a closed, locally symmetric manifold as in Theorem
~\ref{theorem:degree}.  Then $M$ admits no self-maps of degree $>1$.
In particular, $\pi_1(M)$ is {\em co-Hopfian}: every injective
endomorphism of $\pi_1(M)$ is surjective.
\end{corollary}

Note that Corollary~\ref{corollary:hopf} may also be deduced from Margulis's
Superrigidity theorem (for higher rank $M$).  The co-Hopf property for
lattices was first proved by Prasad \cite{Pr}.

\section{The Doaudy-Earle-Besson-Courtois-Gallot map}
\label{section:construction}

\subsection{The construction}
In this section we describe the canonical map $F$ discussed above.  This
construction is due to Duady-Earle
for the hyperbolic plane, was extended by Besson-Courtois-Gallot
\cite{BCG4} to negatively
curved targets, and
extended to symmetric spaces of noncompact type in \cite{CF1,CF2}.
In this section
we will extend this construction further to nonpositively curved
targets with negative Ricci curvature.  For background on nonpositively
curved manifolds and symmetric spaces, see for example \cite{BGS, Eb}.

As above, let $M,N$ be closed, Riemannian $n$-manifolds with $M$
nonpositively curved, and let $f:N\to M$ be any continuous map.  Denote
by $Y$ (resp.\ $X$) the universal cover of $N$ (resp.\ $M$).  Denote by
$\D X$ the visual boundary of $X$; that is, the set of equivalence
classes of geodesic rays in $X$, endowed with the cone topology.  Hence
$X\cup \D X$ is a compactification of $X$ which is homeomorphic to a
closed ball; see, e.g.\ \cite{BGS}.

\bigskip
\noindent
{\bf Idea of the main construction. } Let $\phi$ denote the lift to
universal covers of $f$ with basepoint $p\in Y$ (resp. $f(p)\in X$),
i.e. $\phi=\tilde{f}:Y\to X$.  We first construct a map
$\widetilde{F}:Y\to X$ by ``averaging'' $\phi$ as follows: first embed
$Y$ into a space of measures on $Y$, then push forward each measure via
$\phi$, then smooth out the measure onto $\partial X$ by convolving with
a canonical measure on $\D X$, and finally take the ``barycenter'' of
the resulting measure.  An essential feature of each of these steps is
that they are {\em canonical}, i.e.\ they are equivariant with respect
to the actions of fundamental groups on each of the spaces involved.  It
follows from this that $\widetilde{F}$ descends to a map $F:N\to M$.

Now to define $\widetilde{F}$ precisely.  Actually, it is useful to fix
a parameter $s>h(g)$ for which we define a map $\widetilde{F}_s$.  For a
measure space $Z$, denote by $\meas(Z)$ the space of probability
measures on $Z$.  We denote the Riemannian metric and corresponding
volume form on universal cover $Y$ by $g$ and $dg$ respectively.

Following the method of \cite{BCG1}, we define a
map $\tilde{F}_s:Y\to X$ as a composition $$\xymatrix{ \mathcal{M}( Y)
\ar[r]^{\phi_*} &\mathcal{M}(X)
\ar[r]^{\circledast\nu_z}& \mathcal{M}(\D X) \ar[d]^{\op{bar}} \\ Y
\ar[u]^{\mu^s_x} \ar@{-->}[rr]^{\tilde{F}_s}& & X
\ar@{}[ull]|{\displaystyle{\curvearrowright}} }$$
where the individual maps are defined as follows:

\begin{itemize}
\item The inclusion $Y\to \meas(Y)$, denoted $y\mapsto \mu^s_y$,
is given by localizing the Riemannian volume form $dg$ on $Y$ and
normalizing it to be a probability measure; that is, $\mu^s_y$ is the
probability measure on $Y$ in the Lebesgue
class with density given by
$$\frac{\displaystyle d\mu_y^s}{\displaystyle dg}(z)=\frac{\displaystyle
e^{-sd(y,z)}}{\displaystyle \int_{Y}
  e^{-sd(y,z)}dg}$$
Note that each $\mu_y^s, y\in Y$ is well-defined by the choice of $s$.

\item The map $\phi_\ast$ is the pushforward of measures.

\item The symbol $\circledast\nu_z$ indicates the operation of convolution
with the {\em Patterson-Sullivan measures} $\{\nu_x\}_{x\in X}$
corresponding to $\pi_1(M)<\Isom(X)$ (see \S~\ref{section:ps}).
In other words, the resulting measure $\sigma_y^s$ is defined
on a Borel set $U\subset\partial X$ by
$$\sigma_y^s(U)=\int_{X}\nu_z(U)d(\phi_*\mu_y^s)(z)$$

Since $\Vert \nu_z \Vert=1$, we have
$$\Vert\sigma_y^s\Vert=\Vert\mu_y^s\Vert=1.$$

\item The map $\bary$ is the
{\em barycenter} of the measure $\sigma_y^s$,
which is defined as the unique minimum of a
certain functional on $X$, depending on the measure (see
\S~\ref{section:barycenter}).  The map $\bary$ is
not always defined; when it is defined its domain of definition often
must be restricted.
\end{itemize}

We now describe the last two maps in more detail.  For full
details in the current context we refer the reader to \cite{CF2}.

\subsection{Patterson-Sullivan measures}
\label{section:ps}

Let $\Gamma$ be a discrete group of isometries of a
connected, simply-connected, complete,
nonpositively curved manifold $X$ with no Euclidean direct factors.
Fix a basepoint $p\in X$.

Generalizing work of Patterson and Sullivan, a number of authors,
including Coornaert, Margulis (unpublished), Albuquerque, and Knieper
constructed a remarkable
family $\{\nu_x\}$ of probability measures on $\partial
X$.  These measures, called {\em Patterson-Sullivan
measures}, are meant to encode the density of an orbit
$\Gamma \cdot p$ at infinity as viewed from $y \in {Y}$, giving
the Patterson-Sullivan measure $\nu_y$ on $\D_\infty {Y} $.  See
\cite{Al} and \cite{Kn} for the construction of Patterson-Sullivan
measures in nonpositive curvature.

The key properties of Patterson-Sullivan measures $\nu_x$ are:
\begin{enumerate}
\item {\bf No atoms: }Each $\nu_x$ has no atoms.
\item {\bf Equivariance: }$\nu_{\gamma x}=\gamma_*\nu_x$ for all
$\gamma\in\Gamma$.
\item {\bf Explicit Radon-Nikodym derivative: } For all $x,y\in X$, the
measure $\nu_y$ is absolutely continuous with respect to $\nu_x$.  In
fact the Radon-Nikodym derivative is given explicitly by:

\begin{align}
\label{eqn:R-N}
\frac{d\nu_x}{d\nu_y}( \xi) = e^{h(g)B(x,y, \xi)}
\end{align}

\noindent
where $B(x,y,\xi)$ is the {\em Busemann function} on $X$.  For points
$x,y\in X$ and $\xi\in\D X$, the function $B:X\times X\times\D X
\rightarrow \R$ is defined by
$$B(x,y,\xi)=\lim_{t\rightarrow \infty}d_X(y,\gamma_\xi(t))-t$$
where
$\gamma_\xi$ is the unique geodesic ray with $\gamma(0)=x$ and
$\gamma(\infty)=\xi$.

\item {\bf Support}(Knieper
\cite{Kn}, Albuquerque \cite{Al}): Each $\nu_x$ is
supported on a specific $\Gamma$-orbit in $\D_\infty {X}$.  
\end{enumerate}

The importance of property (3) is that it makes possible
a reasonably explicit computation of the Jacobian of the
map $F_s$; see below.

Property (4) is the most difficult of the properties to prove 
(see \cite{Al} for the locally symmetric case and \cite{Kn} for
the geometric rank one case).  
It's importance is that when
$X$ is a symmetric space of noncompact type, the support of
each $\nu_x$ lies in a specific copy of the {\em Furstenberg boundary}
$\partial_FX$, and thus may be identified algebraically.  This allows
one to convert a geometry problem into an algebra problem; see below.

\subsection{The barycenter functional}
\label{section:barycenter}

For a probability measure on $\partial X$, one can try to define its
{\em barycenter}, or {\em center of mass}, as follows.  One first
defines a map
$$\mathcal{B}:{X}\times \mathcal{M}(\partial X)\to
    \R $$

via
$$\mathcal{B}(x,\lambda) = \int_{\partial_\infty {X}}
    B(x,\theta)d\lambda(\theta)$$
where $B$ is the Busemann function on $X$, based at a fixed basepoint.
Since $X$ has nonpositive sectional curvatures,
the distance function $d$ on $X$ is convex, and hence so is the Busemann
fuction $B$.  It follows that for any
fixed measure $\lambda\in\mathcal{M}(\partial X)$, the
function $\mathcal{B}(\cdot,\lambda)$ is a convex function on $X$.

\begin{proposition}[strict convexity]
\label{proposition:convexity}
Suppose that $X$ has nonpositive sectional 
curvature, negative Ricci curvature uniformly bounded away from $0$, 
and no local $\R$ factors.  
Let $\lambda\in\mathcal{M}(\partial X)$ be a fixed measure.  Suppose
that either 
\begin{itemize}
\item $X$ is a symmetric space of noncompact type, and $\supp(\lambda)$
is the Furstenburg boundary of $X$ 
\item $X$ has geometric rank one, and $supp(\lambda)=\partial X$.
\end{itemize}
Then the function $\mathcal{B}(\cdot,\lambda)$ is strictly convex on $X$.
\end{proposition}

We note that in both cases of Proposition~\ref{proposition:convexity} 
the hypotheses are satisfied when $\lambda$ is Patterson-Sullivan measure.

\begin{proof}
The case when $X$ is a symmetric space of noncompact type is Proposition
3.1 of \cite{CF2}.  Hence we assume that $X$ has geometric rank one 
and negative Ricci curvature (say bounded by $-a^2$).  We also fix a
basepoint $p$ for the Buseann function $B$ on $X$.

As the Busemann $B$ function is convex, the functional 
$$\mathcal{B}(y,\lambda):=\int_{\partial X} B(y,p,\xi)d\lambda(\xi)$$
is convex, being an integral of convex functions.  We must prove that 
$\mathcal{B}(\cdot,\lambda)$ is stictly convex on $X$.

To see this, we first note that while the Hessian of $B(\cdot,p,\xi)$ is only
semi-definite for fixed $\xi$, after taking the integral it becomes
positive definite provided that $D_vd_yB$ does not have a common $0$
direction $w\in S_yY$ for each $v\in S_yY$ which lies in the support of
$\lambda$. Since the
support of $\lambda$ is all of $\partial X$, 
we must show that the tensors $D_vd_yB$ do not
have a common zero $w$. Now we have the decomposition,
$$(D_vd_yB)^*=0\directsum U(v)$$ where $U(v)$ is the second fundamental
form of the horosphere through $v$. In particular the Ricatti equations
imply that $U(v)$ is a semi-definite tensor on the subspace of $T_yY$
orthogonal to $v$. Therefore $w$ must lie in $<v>\directsum \op{ker}
U(v)$ for all $v\in S_yY$, or for all $v\neq
\pm w$ we have $w-<w,v>v \in \op{ker} U(v)$. For such a $w$ whenever
$v\neq w$ then the unique Jacobi field $J(t)$ along the geodesic $v(t)$ with
$v'(0)=v$ satisfying $J(0)=w$ and $J(t)$ bounded for all $t>0$, also
satisfies $J'(0)=0$, i.e. $J(t)$ is parallel along $v$ at $t=0$. In particular the Ricatti equation then implies
that $R(v,w,v,w)=0$. Since this is true for all $v\in S_yY$ with
$v\neq w$ then $R_w:=R(w,\cdot,w,\cdot)$ is the zero tensor, but we assumed
that the Ricci curvature $\op{tr} R_w$ at $w$ was nonzero, and we are
done.
\end{proof}

Note that Example~\ref{example:gluing} shows that
the negativity condition on the Ricci curvature is necessary.

Since we also know that the convex function goes to infinity as $x$ goes
to infinity we can define:

\begin{gather*}
    \op{bar}:\mathcal{M}_1(\partial_\infty {X})\to {X}\\
    \op{bar}(\lambda):=\text{unique minimum of }
    \mathcal{B}(\cdot,\lambda) \\
\end{gather*}

\begin{figure}[ht]
\begin{center}
\psfrag{b}{\hspace{.7in} bar($\mu$)}
\psfrag{n}{$\mu$ on $\partial X$}
\scalebox{1}{\includegraphics{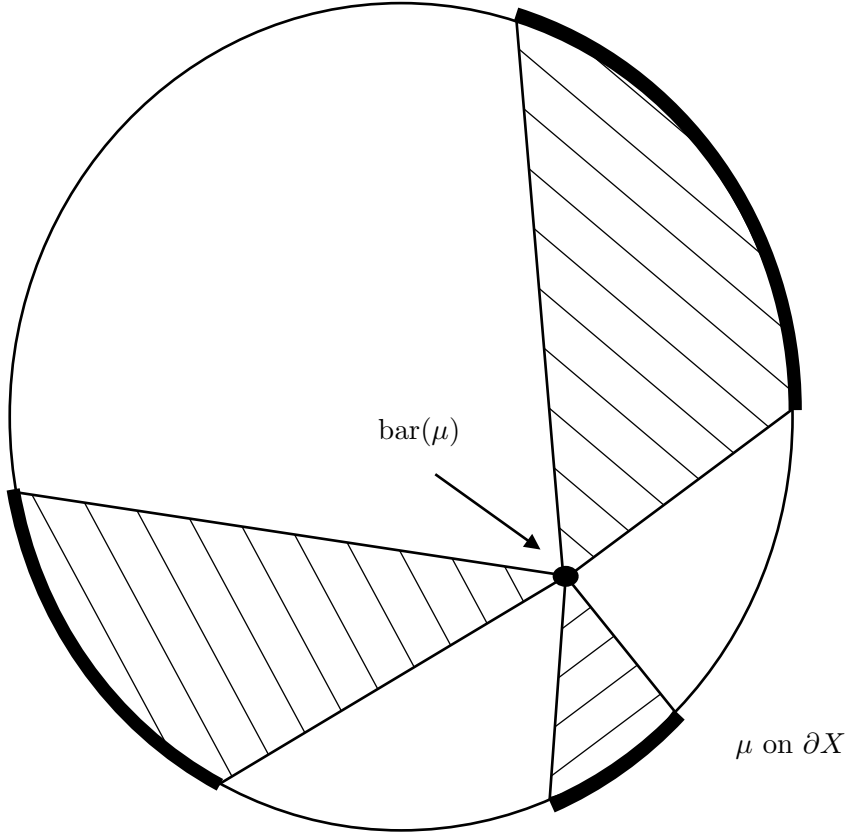}}
\end{center}
\caption{The barycenter of a meaure.}
\end{figure}

We can now define $\tilde{F}:Y\to X$ by
   $$ \tilde{F}(x):=\op{bar}(\tilde{f}_\ast (\mu_x^s)\circledast
\nu_x)$$
as described above.  From the equivariance properties described above,
it follows that $\tilde{F}$ is equivariant with respect to the
homomorphism $f_\ast:\pi_1(N)\to \pi_1(M)$, so that $\tilde{F}$ descends
to a map $F:N\to M$.

\section{Bounding the Jacobian}
\label{section:jacobian:bound}

The power of the map $F$ comes from the fact that one can often obtain
explicit estimates on its Jacobian.  As we saw in
\S~\ref{section:degree:argument}, such estimates are the key to all
applications of the barycenter method.  In this section, we sketch how
estimates on $|\Jac(F)|$ have been obtained in some cases.

\subsection{Existence of some bound}

Let $\sigma_y^s=\tilde{f}_\ast (\mu_x^s)\circledast
\nu_x$.  Now $F_s$ is defined by the
implicit vector equation:
\begin{equation}
\label{eq:implicit}
\int_{\D_\infty {X}}
 d{B}_{(F_s(y),\theta)}(\cdot) d\sigma_y^s(\theta)=0
\end{equation}

Differentiating (~\ref{eq:implicit})
and applying the Implicit Function Theorem then gives a
formula for $|\Jac F_s|$.  An application of H\"{o}lder's inequality and
some further estimates then give:

\begin{equation}
\label{eq:Jacobian}
| \Jac F_s|\leq
\left(\frac{s}{\sqrt{n}}\right)^n\frac{\det\left(\int_{\D_\infty X}
    d{B}_{(F_s(y),\theta)}^2(\cdot )
    d\sigma^s_y(\theta)\right)^{1/2}}{\det\left(\int_{\D_\infty X}
    Dd{B}_{(F_s(y),\theta)}(\cdot,\cdot)
    d\sigma^s_y(\theta)\right)}
\end{equation}

Note that all of this works as long as $F_s$ is well-defined.

We wish to bound the right hand side of (~\ref{eq:Jacobian}).  In order to do
this we first convert the problem to a Lie groups problem.  This is
where we use crucially the assumption that $M$ is locally symmetric.
An easy argument allows us to assume that the universal cover of $M$ is
irreducible, which we now do.

Recall Property (4) of Patterson-Sullivan measures, which states that in
the locally symmetric case, the support of each Patterson-Sullivan is
contained in a specific $\Gamma$-orbit which may also be identified with
the {\em Furstenberg boundary} corresponding to the Lie group
$G=\Isom(\tilde{M})$ of isometries of the symmetric space $\tilde{M}$.
From this one can deduce the following:

\bigskip
\noindent
{\bf Key fact: }Can replace $\partial_\infty X$ in (~\ref{eq:Jacobian})
by the maximal compact subgroup $K\subset G=\Isom(X)$.
\bigskip

To further simplify the right hand side of (~\ref{eq:Jacobian}), we combine :
\begin{itemize}
\item Eigenvalue estimates on $DdB$, and
\item The fact that, for $M_i$ positive semidefinite,
$\det(\sum M_i)$ is a nondecreasing homogeneous polynomial in
eigenvalues of the $M_i$.
\end{itemize}
to obtain
\begin{equation}
\label{eq:lastequation}
|\Jac F| \leq C\frac{\left(\det \ds{\int}_{K}
    O_{\theta}\begin{pmatrix} 1 & 0 \\
      0 & 0 \end{pmatrix}O_\theta^* \
    d\sigma_y(\theta)\right)^{\frac12}}{\det\ \ds{\int}_{K}
  O_{\theta}\begin{pmatrix} 0 & 0 \\ 0 & I
   \end{pmatrix}O_\theta^* \ d\sigma_y(\theta)}
\end{equation}
where $I$ is the identity matrix of size $n-\rank(X)$, the matrix $O_\theta$
is an element of $K$, and the constant $C$ depends only on $n$ and on
the Ricci curvatures of $M,N$.  We remark that the eigenvalues of $DdB$
are determined purely by algebraic data attached to the Lie algebra of $G$.

The strategy now is to bound the right hand side of
(~\ref{eq:lastequation}) {\em over all probability measures} on $K$.  In
some sense this is the weakness of the barycenter method: without using
additional information, one does not always obtain a bound as needed to
prove Theorem~\ref{theorem:degree}, much less 
Conjecture~\ref{conjecture:entrig}, at least for certain
locally symmetric manifolds; such examples which show that the method
can fail are given in \S~\ref{section:counterexamples} below.

Some bounds on the right hand side of (~\ref{eq:lastequation})
are possible in most cases, however; the main technique is
called ``eigenvalue matching'' in \cite{CF2}.  The idea is that
for each small eigenvalue
in the denominator, one tries to find two comparably small eigenvalues in the
numerator.  The hard part is to make ``comparably small'' independent of the
measure.  This involves a detailed analysis of action
of the maximal compact subgroup $K$ on
subspaces of Lie algebra of $G=\Isom(X)$ (a simple Lie
group). 

The first basic fact used here is that the kernel of the operator
in the numerator of the right hand side of
(~\ref{eq:lastequation}) contains the cokernel of the operator in
the denominator.  The second main ingredient is the following
proposition, which gives a $K$-invariant subspace perpindicular
to $K$-orbit of a flat, and of twice dimension of flat
(=$\rank(X)$). This fact is precisely where the hypothesis that
$G\neq \SL(3,\R)$ is used.

\begin{proposition}[Eigenspace matching]
There is a constant $C$, depending only on $\dim
{X}$ such that for any subspace $V\subset T_{x}{X}$ with $\dim
V\leq\rank({X})$ there is a subspace $V'\subset V^\perp$ of dimension
$2\cdot \dim V$ satisfying
$$\measuredangle (O_\theta\cdot
V',A^\perp) \leq C\measuredangle (O_{\theta}\cdot V,A)$$
for all $O_\theta\in K$.
\end{proposition}

\begin{figure}[ht]
\begin{center}
\scalebox{1}{\includegraphics{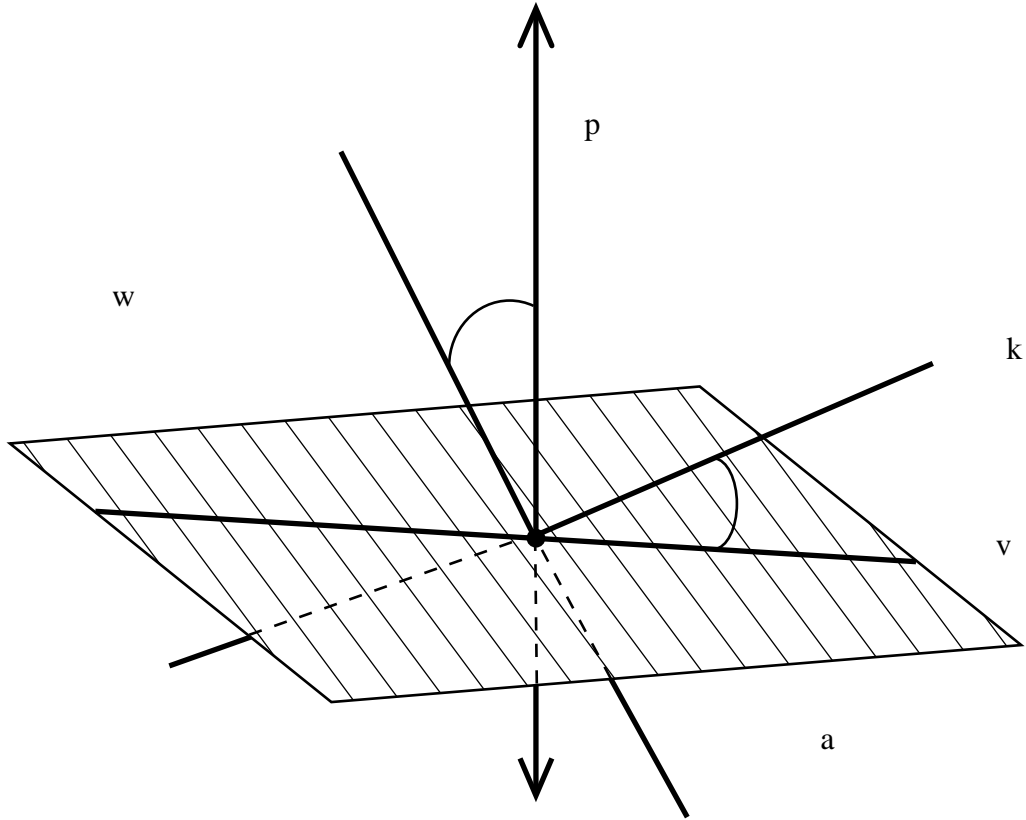}}
\end{center}
\caption{Eigenvalue matching.}
\end{figure}

\subsection{The perfect bound in rank one}

When $M$ has real rank one, i.e.\ when $M$ is negatively curved, 
Besson-Courtois-Gallot are able to obtain the infinitessimal
rigidity described on page \pageref{equation:infinitessimal}
(and hence a proof of the entropy rigidity conjecture) because of
the following remarkable phenomenon.  It is possible to write the
operator in the denominator of the right hand side of the inequality
(~\ref{eq:lastequation}) explicitly in terms of the operator in
the numerator and the complex structure on the symmetric space;
see \cite{BCG2} for details.  For example, for the case when $M$
is real hyperbolic, the denominator operator is the identity
minus the numerator.  They are then able to solve the
minimization problem explicitly by using (quite nontrivial)
linear algebra.

\section{Negative results in various settings}
\label{section:negative}

In this section we describe in greater detail some negative results in
the settings given in the table on page \pageref{table:settings}

\subsection{Finsler metrics}

Given a manifold, $M$, if we equip each tangent space $T_xM$ with a norm
$F_x$ which depends smoothly on $x$ then this gives rise to a Finsler
metric by taking the distance between two points $x$ and $y$ to be the
infimum of lengths of curves connecting $x$ to $y$ where length is
measured using the norm on their tangent vectors. For our purposes we
shall assume that the unit ball in each $T_xM$ given by the norm is
strictly convex with $C^2$ boundary. This is a common assumption since it
allows us to form a positive definite inner product, $g_u$, for each $u\in
T_xM$ whose matrix in local coordinates $(x_i,\dot{x}_i)$ for $TM$ is
given by entries $$\frac{\partial^2 F^2}{\partial x_i\partial x_j}(u).$$
As expected, when $F$ comes from a Riemannian inner product then the $g_u$
all equal the inner product on $T_xM$. Geodesic flow makes sense for $F$
as it does for the Riemannian case, and therefore have a dynamical notion
of topological entropy.

Let $B_x^F(R)$ be the ball of radius $R$ in $T_xM$ in the norm $F$. For any
Riemannian metric $g$ on $TM$ we can similarly define $B_x^g(R)$. Then the
volume form for $F$ is defined at each point $x$ to be,
$$dF(x):=\frac{\Vol_g B_x^g(1)}{\Vol_g B_x^F(1)}dg(x),$$
which is independent of the choice of inner product $g$. In particular, the
quantities $\Vol_F(M)$ and $h_{\vol}(F)$ make sense.

For deformations of a single Finsler metric, we can keep $dF$ fixed by
simply keeping the $\Vol_g B_x^F(1)$ fixed for each $x$ relative to a
fixed inner product $g$. Keeping the topological entropy constant as
well is more subtle.  P. V\'{e}rovic (\cite{Ve}) has shown for any
volume preservering deformation $F_t$ for $t\in (-\eps,\eps)$ of a
compact hyperbolic manifold through Finsler metrics where where $F_0$ is
a hyperbolic Riemannian metric,
$$\left.\frac{d}{dt}\right|_{t=0}\op{h_{top}}(F_t)=0.$$
For surfaces
he shows the same result for deformations which preserve the Liouville
volume of each unit tangent sphere.  Moreover, there exist
one-parameter families $F_t$ where topological entropy can be held
constant; in particular the rigidity part of Theorem~\ref{theorem:bcg}
does not hold when the class of metrics is extended to Finsler
metrics.

In higher rank things get worse: V\'{e}rovic also proved in \cite{Ve}
that for compact
higher rank symmetric spaces $(M,g_0)$, there exists a Finsler metric $F_0$
(invariant under all local isometries) with $\Vol_{F_0}(M)=\Vol_{g_0}(M)$ and
$h_{\vol}(F_0)<h_{\vol}(g_0)$; that is, the locally symmetric metric
on a higher rank manifold does {\em not} minimize normalized
entropy among Finsler metrics.

Any invariant Finsler metric is determined by its values on a single
maximal flat $A$. In fact, if $\mathfrak{a}$ is the corresponding maximal
abelian subalgebra and $\Lambda^+$ is a set of positive roots $\alpha$
each with multiplicity $m_\alpha$ then $F_0$ is given by

$$F_0(H):=\sum_{\alpha\in\Lambda^+_\mathfrak{a}}m_\alpha |\alpha(H)|$$
for each $H\in \mathfrak{a}$.

In \cite{BCG3} Besson, Courtois and Gallot asked the following
intriguing question.
\begin{question}
Does V\'{e}rovic's Finsler metric $F_0$ minimize normalized entropy
among all Finsler metrics?
\end{question}

\subsection{Magnetic field flows}
Let $(M,g)$ be a Riemannian manifold, let
$\pi:TM\to M$ be natural projection from the tangent bundle $TM$ of
$M$.  The metric $g$ determines an isomorphism of $TM$ with the
cotangent bundle; let $\omega_0$ be the symplectic 2-form on $TM$
formed by pulling back the
canonical symplectic form on $T^*M$ by this isomorphism.
Now consider any closed 2-form
$\Omega$ on $M$. Then for any $\lambda\geq 0$ it is easy to verify that
$$\omega_\lambda=\omega_0+\lambda
\pi^*\Omega$$
is also a symplectic form on $TM$.

Consider the usual Hamiltonian
$$H_x(v):=\frac12 g_x(v,v)$$
and define the {\em $\lambda$-magnetic flow}\/ of the pair $(g,\Omega)$, denoted
$\phi_\lambda$,  to be the
Hamiltonian flow of $H$ with respect to $\omega_\lambda$. Consider the bundle
map $Y:TM \to TM$ given implicitly by
$$\lambda\Omega_x(u,v)=g_x(Y_xu,v)$$
Then for $v\in TM$ the curve $t\mapsto
\phi^t_\lambda(v)=(\gamma(t),\dot{\gamma}(t))$ is characterized by the fact
that
$$\nabla_{\dot\gamma(t)}\dot\gamma(t)=Y_{\gamma(t)}\dot\gamma(t),$$
where $\nabla$ is the connection for the metric $g$; in
particular $\phi_0^t$ is just the geodesic flow for $g$.

Now we consider the case where $M$ admits an Anosov geodesic flow.
For instance we could take $M$ to be any negatively curved compact
manifold. In \cite{PaPa1} and \cite{PaPa2}, G. and M. Paternain show
that for $\Omega\neq 0$ then $\htop(\phi_\lambda)$ is strictly
decreasing for $\lambda\geq 0$ in the (nonempty) interval containing
$\lambda=0$ for which $\phi_\lambda$ is an Anosov flow.  K. Burns and
G. P. Paternain further show that $\htop(\phi_\lambda)$ is decreasing
between "Anosov intervals" (see \cite{BP}). In this paper the authors
also give an example of a higher genus surface $M$ with a certain
choice of $\Omega$ where $\htop(\phi_\lambda)$ increases for $\lambda$
in between disjoint intersvals where $\phi_\lambda$ is Anosov.
Therefore, while metrics corresponding to local maxima and minima for
$\htop(\phi_\lambda)$ may be speciual, no standard entropy rigidity
exists for such flows.

\subsection{Metric entropy}

Lastly, we mention the work of L. Flaminio in \cite{Flaminio}. The
{\em metric entropy} is the measure theoretic entropy of the Liouville
measure on $SM$. The measure theoretic entropy of the Bowen-Margulis
measure is the same as the topological entropy, which in turn is the
same as the volume entropy. Hence for hyperbolic metrics, for instance,
the metric entropy and volume entropy coincide. In general the
topological entropy is always greater than or equal to the metric
entropy. Flaminio shows that the metric entropy for volume
preserving deformations of a hyperbolic metric on a three manifold $M$
does not have a local maximum at the hyperbolic metric. This is in
contrast to the two dimensional case where A. Katok \cite{Katok}
showed that the metric entropy does have a local maximum there.

\section{Positive results in various settings}
\label{section:positive}

In this section we describe some known cases where a locally symmetric
metric structure minimizes $h(g)$ for a suitable normalization.

\subsection{Finsler metrics}

Recall the definition of a Finsler metric given above. We call a Finsler
metric $F$ {\em reversible} if it satisfies $F(-v)=F(v)$ for all $v\in TM$.
Moreover we define an {\em eccentricity factor} $N(F)$ for a Finsler metric $F$ on a
manifold $X$ to be
$$N(F):=\max_{x\in X}\max_{u\in S_x^F(1)}\max_{v\in S_x^{g_u}(1)}
 \frac{F_x(v)^n\Vol_{g_u}(B_x^F(1))}{\Vol_{g_u}(B_x^{g_u}(1))},$$
where $S_x^F(1)$ and $S_x^{g_u}(1$ are the unit spheres in the norms $F$ and
$g_u$ respectively. For this setting we redefine the volume entropy functional
$\ent$ to be
$$\ent(F):=N(F)h(F)^n\vol_F(X)$$
where $n$ is the dimension of $X$.

To each direction in $T_xX$ we may ascribe a Riemannian metric,
and consider the corresponding curvature tensor. These are the
flag curvatures (see \cite{BaoChern}).  With this notation, J.
Boland and F. Newberger \cite{BolandNewberger} proved the
following theorem.

\begin{theorem}[Finsler entropy rigidity]
Let $(X,g_0)$ be a compact, n-dimensional, locally symmetric Riemannian
manifold of negative curvature $(n\geq 3)$ and $(Y,F)$ a compact, reversible,
Finsler manifold of negative flag curvature that is homotopy equivalent to
$(X,g_0)$. Then
$$1\leq \frac{\ent(F)}{\ent(g_0)}$$
with equality if and only if $F$ is homothetic to $g_0$
\end{theorem}

Notice that if $F$ is Riemannian then $N(F)=N(g_0)=1$ and so this reduces
to the usual entropy rigidity theorem for the case of maps homotopic to the
identity.

\subsection{Foliations}

Let $N$ and $M$ be compact topological manifolds supporting continuous
foliations $\mathcal{F}_{N}$ and $\mathcal{F}_{M}$ by leaves which are smooth
Riemannian manifolds, and such that the metrics on the leaves vary
continuously in the transverse direction. We suppose that the leaves of
$\mathcal{F}_{M}$ are locally isometric to $n$-dimensional symmetric spaces of
negative curvature, $n \geq 3$; by continuity of the metrics these are all
locally homothetic to a fixed symmetric space $(\tilde{X}, g_{o})$.

For any leaf $(L,g_L)$, we may define
$$\bar{h}(g_L)=\limsup_{R\to\infty} \frac{\Vol(B(x,R))}{R}$$
and similarly we can define $\underline{h}(g_L)$ as the
$\liminf_{R\to\infty}$.  These numbers do not depend on the choice of
$x$.  We then define the {\em volume growth entropy} $h(g_{L})$ to be
$$h(g_{L}) = \inf\left\{s>0
  \left| \int_0^\infty e^{-s t} \Vol S(x,t) dt <\infty \right.
\right\},$$ where $S(x,t)$ is the sphere of radius $t$ about $x$ in the
universal cover $\tilde{L}$ of $L$.  This quantity is independent of $x
\in L$.  Hence we can define a function $f:N\to [0,\infty]$ by letting
$f(x)$ be $h(g_L)$ for the leaf $L$ containing $x$.  Of course $f$ is
contant on each leaf.

The function $f$ is also measurable.  This follows from the fact
that the transverse
continuity of the leafwise metrics implies that for each $R$, the
function
$$x\mapsto \int_0^{R} e^{-s t} \Vol S(x,t) dt$$
is continuous on $N$. On
$(M,\mathcal{F}_{M})$ the entropy is constant and we denote it by $h(g_{o})$.
From the definition of $h(g_{L})$ it follows almost immediately that
$$\underline{h}(g_L)\leq h(g_{L}) \leq \bar{h}(g_L).$$

For the foliation $(N,\mathcal{F}_{N})$ we assume that the leaves $(L, g_{L})$
are strictly negatively curved, and which satisfy the stronger
condition that for all $x,y\in L$ there is a $\delta <h(g_{L})$ such that
\begin{equation}\label{eq:vol-condition}
\limsup_{R\to\infty} \frac{\Vol S(x,R)}{\Vol S(y,R)} \leq C e^{\delta d(x,y)},
\end{equation}
where $S(x,R)$ is the sphere of radius $R$ in $L$. Such a leaf will be called
a {\em Patterson-Sullivan manifold}.

Finally, assume we have a leaf-preserving homeomorphism
$$f:(N, \mathcal{F}_{N}) \rightarrow (M,\mathcal{F}_{M})$$
which is leafwise $C^1$ with transversally continuous leafwise derivatives,
but not necessarily transversally differentiable.

Equip the foliation $(N, \mathcal{F}_{N})$ with any choice of finite,
transverse, holonomy quasi-invariant measure $\nu$ (see Hurder \cite{Hurder}
or Zimmer \cite{Zimmer} for the definition and existence). Holonomy
quasi-invariance simply means that the push forward of $\nu$ under any
holonomy map is in the same measure class as $\nu$. This measure $\nu$
provides us with a global finite measure $\mu_{N}$ on $N$ which is locally a
product of $\nu$ with the Riemannian volumes $\op{dvol}_{L}$ of the leaves
$L$.

\begin{theorem}[Foliated entropy rigidity, I]
\label{thm:main3}
Let $(N, \mathcal{F}_{N})$ be a continuous foliation of the compact manifold
$N$ such that $\nu$-almost every leaf is a Patterson-Sullivan manifold.
Suppose that $f: (N, \mathcal{F}_{N}) \to (M,\mathcal{F}_{M})$ is a
foliation-preserving homeomorphism, leafwise $C^1$ with transversally
continuous leafwise derivatives, and that $f_*\nu$-almost every leaf of
$(M,\mathcal{F}_{M})$ is a rank one locally symmetric space. Then there exists
a finite measure $\mu_{M}$ on $M$ which is locally the product of
$\op{dvol}_{o}$ with a transverse quasi-invariant measure $\nu_{o}$ such that
$$\int_{M} h(g_{o})^n d\mu_{M} \leq \int_{N} h(g_{L})^n
d\mu_{N},$$ and equality holds if and only if $\nu$-almost every leaf
$(L,g_{L})$ is homothetic to its image $(f(L), g_{o})$.
\end{theorem}

When the foliation admits a holonomy invariant measure $\nu$ then we may take
$\nu_o=f_*\nu$. When $\nu$ is just holonomy quasi-invariant however, then
$\nu_o$ is the push forward of $\nu$ under the natural map $F$ defined below.

When the foliation $(N, \mathcal{F}_{N})$ is ergodic with respect to $\nu$,
then the entropy function $h(g_{L})=h(g)$ is constant on $N$, and we get the
following.
\begin{corollary}[Foliated entropy rigidity, II]
\label{cor:ergodic}
  Under the same assumptions as in the main theorem, if $(N,
  \mathcal{F}_{N})$ is ergodic, then $h(g_{o})^n \Vol(M, \mu_{M}) \leq
  h(g)^{n} \Vol(N, \mu_{N})$ with equality if and only if $\nu$-almost
  every leaf $(L,g_{L})$ is homothetic to $(f(L), g_{o})$.
\end{corollary}

\bigskip
\noindent
{\bf Remark.} If $(N, \mathcal{F}_{N})$ and $(M,\mathcal{F}_{M})$ are
foliations such that almost every leaf is compact or simply connected,
then the requirement that the homeomorphism $f$ be leafwise $C^1$ can be
dropped.  In particular if the foliations have just one leaf and $\dim N
\ne 3,4$, any homotopy equivalence induces a homeomorphism between $N$
and $M$ (see \cite{FarrellJones}). Therefore, when $\dim N \ne 3,4$,
Corollary~\ref{cor:ergodic} recovers Theorem~\ref{theorem:bcg}.

\bigskip
Of course one can also ask for foliated versions of Conjecture
~\ref{conjecture:entrig} as well.

\subsection{Finite volume manifolds}

For the case of finite volume manifolds J. Boland, J. Souto and the second
author showed the following analogue of the Real Schwarz Lemma.

\begin{theorem} [Volume Theorem \cite{BCS}]
\label{thm:bcs}
Let $(M,g)$ and $(M_o,g_o)$ be two oriented complete finite volume riemannian
manifolds of the same dimension $n \ge 3$ and suppose that
$$Ric_g\ge -(n-1)g,\ \ \ \hbox{and}\ \ \ -a\le K_{g_o}\le -1.$$
Then for all proper continuous maps $f: M\to M_o$,
$$\vol(M,g)\ge \vert\deg(f)\vert\vol(M_o,g_o),$$
and equality holds if
and only if $M$ and $M_o$ are hyperbolic and $f$ is proper homotopic
to a Riemannian covering.
\end{theorem}

When $M$ and $M_o$ are compact, Theorem~\ref{thm:bcs} follows from a
real Schwarz lemma proved by Besson, Courtois, and Gallot in
\cite{BCG5}. As in the case of all noncompact situations, to apply the
method of \cite{BCG1} and
\cite{BCG5}, the fundamental difficulty is proving the properness of the
natural map.

Theorem~\ref{thm:bcs} implies the following.

\begin{corollary}
Under the hypotheses of Theorem~\ref{thm:bcs},
$$\op{Minvol}(M)\ge \deg(f) \vol(M_o).$$
\end{corollary}

Now restrict $f$ to be degree 1, and consider what happens when
$\op{Minvol}(M)=\vol(M_o).$ Bessieres proved in \cite{Be1} that if
there is a degree 1 map $f: M\to M_o$ from the compact $n$-dimensional
$M$ to the compact $n$-dimensional real hyperbolic $M_o$ such that
$\op{Minvol}(M)=\vol(M_o)$, then $M$ and $M_o$ are diffeomorphic. He
also gave examples in \cite{Be2} of a finite volume manifold $M$ with
the same simplicial volume as a hyperbolic manifold $M_o$ such that
there is a degree 1 map from $M$ to $M_o$ and $\op{Minvol}(M) \leq
\vol(M_o)$ but $M$ and $M_o$ are not even homeomorphic. By the above
corollary, $\op{Minvol}(M) =\vol(M_o)$. However, given such an example,
our next result shows that the pointed Lipschitz limit of some
subsequence of any sequence of metrics whose volumes achieve
$\op{Minvol}(M)$ is isometric to $M_o$.

\begin{theorem} [Hyperbolic \boldmath$\Minvol$ rigidity, finite volume case]
\label{optimal}
Let $M$ and $M_o$ be finite volume manifolds of the same dimension $n \geq 3$,
$M_o$ real hyperbolic, and $f: M \to M_o$ a continuous, proper, degree 1 map.
If $\op{Minvol}(M)= \vol(M_o)$, then for any sequence of metrics $g_i$
realizing the minimal volume of $M$, there are $p_i\in M$ and a subsequence
$g_{i_j}$ such that $(M,p_{i_j},g_{i_j})$ converges in the pointed Lipschitz
topology to Riemannian manifold isometric to $M_o$.
\end{theorem}

In particular, the topology of the limit manifold changes to that of
$M_o$.
 
The last main result of \cite{BCS} is the finite volume version of the entropy
rigidity result found in \cite{BCG1}.

\begin{theorem}[Rank one entropy rigidity, finite volume case]
\label{thm:entropy3}
Let $(M,g)$ be an $n$-dimensional finite volume manifold of nonpositive
sectional curvature, $n \geq 3$, and $h(g)$ its volume growth entropy. Let
$(M_o,g_o)$ be an $n$-dimensional finite volume rank one locally symmetric
manifold and $h(g_o)$ its volume growth entropy. If $f:M \to M_o$ is a
continuous, proper map of degree $\deg(f)>0$, then
$$h(g)^n \vol(M,g) \geq
\deg(f)h(g_o)^n \vol(M_o,g_o)$$
and equality holds if and only if $f$ is proper
homotopic to a Riemannian covering.
\end{theorem}

As in \cite{BCG2}, this gives a quick proof of the Mostow rigidity theorem
for finite volume negatively curved locally symmetric manifolds.

\subsection{Quasifuchsian representations}

The barycenter method has recently been applied by
Beson-Courtois-Gallot \cite{BCG5} to convex cocompact
(infinite covolume)  representations, generalizing Bowen's
rigidity theoem about quasifuchsian groups.  We now
describe one of the several results in \cite{BCG5} along
these lines.

Let $X$ be a negatively curved manifold.  A faithful
representation $\rho:\Gamma\to \Isom(X)$ is {\em convex
cocompact} if $\rho(\Gamma)$ acts cocompactly on the
convex hull of its limit set $\Lambda(\rho(\Gamma))$.  
Attached to each $\rho$ is a number, namely the Hausdorff
dimension $\dim_{\mathcal{H}}(\rho(\Lambda(\Gamma)))$.  
Amazingly, this single number can be used to characterize
the {\em totally geodesic} representations, i.e.\ those
representations leaving invariant a totally geodesic
submanifold.  This theorem was originally proved by Pansu,
Bourdon, and Yue, generalizing an earlier theorem of
Bowen.  Besson-Courtois-Gallot gave another proof in
\cite{BCG5} using the barycenter method.

\begin{theorem}[Quasifuchsian rigidity]
\label{theorem:quasifuchsian} Let $\Gamma$ be the
fundamental group of a closed, hyperbolic $n$-manifold,
and let $X$ be a connected, simply-connected manifold with
sectional curvature $\leq -1$.  Then
$$\dim_{\mathcal{H}}(\rho(\Lambda(\Gamma)))\geq n-1$$ with
equality if and only if $\rho$ is totally geodesic.
\end{theorem}

There is also a version for complex hyperbolic manifolds
(see cite{BCG5}).  The idea of the proof is to define the
canonical map as for maps between $n$-manifolds, and to
bound the Jacobian of this map; Theorem
~\ref{theorem:quasifuchsian} then follows reasonably
quickly.  The main difficulty is that the dimensions of
the domain and target spaces are different.  To repair
this one notes tha $\rho$ gives a quasi-isometric
embedding between negatily curved spaces, hence a
homeomorphic embedding of one boundary at infinity into
another.  The general technique can then be carried out
once one finds the right notion of volume.

\subsection{Alexandroff Spaces}

In this section we report on recent work of P. Storm, who extended
some of the results of \cite{BCG1} to certain Alexandroff space
domains within the same bilipschitz class as a fixed hyperbolic
manifold.  All of the results can be found in \cite{St}. 

In brief, an Alexandroff space $X$ with curvature bounded below by $-1$ is
any complete locally compact metric space with finite Hausdorff
dimension such that every point $x\in X$ has a neighborhood containing
a geodesic triangle such that the comparison triangle in $\hyp^2$ with
the same side lengths has the property that the distance form any
vertex of the triangle to any point on the opposite side is shorter
than the corresponding distance on the triangle in $X$. Storm's main
result is the following.

\begin{theorem}[Alexandroff domain]
\label{thm:storm} Let $(X,d)$ be an Alexandroff space with curvature
bounded below by $-1$, and $Y$ a closed, hyperbolic $n$-manifold, $n\geq
3$. If $X$ and $Y$ are bilipschitz, then
$$\Vol(X)\geq \Vol(Y)$$
\end{theorem}

The proof of this theorem is at first glance unrelated to the method
of Besson, Courtois, and Gallot as outlined in the opening sections of
this paper. Nevertheless, it is actually of the same vein. The main
technique is to extend the idea of \emph{Spherical Volume} as
presented in \cite{BCG1} to this new setting. Instead of going through
the (perhaps insurmountable) difficulty of forming a barycenter map
from $X$ to $Y$, one can in a sense stop ``half-way'' and attempt to
do the analysis in $\mathcal{M}_1(X)$ and $\mathcal{M}_1(X)$. Since
this embedding is in essence an $L^1$ approach it is preferrable to
modify this to an $L^2$ approach which leaves us working in a Hilbert
space.

More precisely, one embeds $\tilde{X}$ into $L_2(\tilde{X})$ in a
$\pi_1$ equivariant way via the map $\Phi_s(x):=e^{s/2 d(x,\cdot)}$
for any $s>h(X)$, and then projects this to the unit Hilbert sphere
$S^\infty(\tilde{X})$.  We recall that this is the embedding achieved
by using the square root of the Radon-Nikodym derivatives of the
unormalized $\sigma_y^s$, which are in $L^2$ by definition of the
volume growth entropy $h(X)$. If $f:\tilde{Y}\to \tilde{X}$ represents
a bilipschitz $\pi_1$-equivariant map given by the hypotheses of the
theorem, then the map $\tilde{I}:L_2(\tilde{X})\to L_2(\tilde{Y})$
given by $\mathcal{I}(g):=(g\circ f)\Jac(f)$ is a $\pi_1$ equivariant
isometry.  Hence it restricts to $\mathcal{I}:S^\infty{\tilde{X}}\to
S^\infty{\tilde{Y}}$. The composition $\tilde{F}:=\mathcal{I}\circ
\Phi_s\circ f:\tilde{Y}\to S^\infty{\tilde{Y}}$ descends to a map
$F:Y\to S^\infty(Y)$ with image in the positive orthant. 

For such maps
there is a notion of volume, and without defining it, we simply point out that
spherica volume is the infimum of the volume of such maps. In
particular 
$$\vol(F)\geq \op{Spherical}\vol(Y)=\left(\frac{(n-1)^2}{4n}\right)^{n/2}\Vol(Y).$$
The last equality was proven by Besson Courtois and Gallot in \cite{BCG1}.

Since Alexandroff spaces are in a sense Riemannian manifolds off of a
measure zero set, and $Y$ is a Riemannian manifold, we may use
standard integration theory. Of course we now come to the difficult
(and omitted) part of the proof of Theorem~\ref{thm:storm} which is to
show
$$\vol(F)\leq \left(\frac{s^2}{4 n}\right)^{n/2}\vol(X).$$
It then follows that,
$$\vol(X)\geq \left(\frac{n-1}{s}\right)^n\vol(Y).$$
Lastly, Burago,
Gromov and Perel'man prove that for Alexandroff spaces with curvature
bounded below by $-1$, $h(X)\leq n-1$ and hence Theorem
~\ref{thm:storm} follows.

The main applications of Theorem~\ref{thm:storm} are to two
important classes of Alexandroff spaces which arise naturally in the study of
hyperbolic manifolds: doublings of convex cores and cone manifolds.

\subsubsection{Convex Cores}
The (metric) double $DC_M$ of the convex core
$C_M$ of a convex cocompact manifold $M$ is the
result of identifying the boundaries of two copies of $C_M$ and then
extending the induced metric on each copy of $C_M$ to the whole. While
as a topological manifold $DC_M$ can always be smoothed (e.g. by
taking the double a neighborhood of $C_M$), the point is that any
Riemannian metric on the resulting smooth manifold cannot agree with
the metric on each copy of $C_M$ as a subspace, unless these had
totally geodesic boundaries to begin with.  On the other hand, $DC_M$,
naturally carries the structure of an Alexandroff space with
Alexandroff curvature bounded below by $-1$ and the above theorem
applies.

Let $\text{CC} (N)$ be the space of complete convex cocompact hyperbolic
manifolds diffeomorphic to the \textit{interior} of smooth compact
$n$-manifold $N$. In analogy with the volume entropy results discussed
in this paper, one wants to know that the topological invariant
$$\mathcal{V} (N) := \inf_{M \in \text{CC}(N)} \{ \text{volume of the convex core } C_M \text{ of } M \},$$
is minimized (uniquely?) for certain special hyperbolic
manifolds.

It is consequence of Thurston's Geometrization and Mostow Rigidity
that for $n=3$, $N$ is acylindrical if and only if there exists unique
$M_0 \in \text{CC}(N)$ such that $\partial C_{M_0}$ is totally
geodesic \cite{Th2}. Combining this with Theorem~\ref{thm:storm},
Storm obtains the following.

\begin{corollary}
\label{3-cor}.   Let $N$ be an acylindrical compact  irreducible $3$-manifold
such that $\text{CC} (N)$ is nonempty.  Then there exists a unique
$M_0 \in \text{CC}(N)$ such that $\mathcal{V} (N) =
\text{Vol}(C_{M_0})$. Moreover, $\partial C_{M_0}$ is totally
geodesic.
\end{corollary}

In fact he shows that any $C_{M_0}$ with totally geodesic boundary is
the unique minimizer of $\vol(C_N)$ among hyperbolic $N$ homotopy
equivalent to $M_0$. (Note we are using that $M_0$ is tame). We should
also mention that Bonahon had previously shown in \cite{Bon} that
$M_0$ is a strict local minimum of $\vol(C_N)$.

\begin{remark}
Storm is able to also obtain results giving exact relations between the
Gromov norm of such $DC_{N}$ and covers in terms of  $\mathcal{V}(N)$.
\end{remark}

\subsubsection{Application to Cone Manifolds}\label{subsec:cone}

Theorem~\ref{thm:storm} may also be applied to cone-manifolds with all
cone angles $\le 2\pi$.

\begin{Definition}
\cite[pg.53]{CHK} An \emph{$n$-dimensional cone-manifold} is a topological
manifold, $M$, which admits a triangulation giving it the structure of a PL
manifold (i.e. the link of each simplex is piecewise linear homeomorphic to a
standard sphere) and $M$ is equipped with a complete path metric such that the
restriction of the metric to each simplex is isometric to a geodesic simplex
of constant curvature $K$. The singular locus $\Sigma$ consists of the points
with no neighborhood isometric to a ball in a Riemannian manifold.
\end{Definition}

It follows that

{\noindent}$\bullet \ \  \Sigma$ is a union of totally geodesic closed
simplices of dimension $n-2$.

{\noindent}$\bullet \ $   At each point of $\Sigma$ in an open
$(n-2)$-simplex, there is a \emph{cone angle} which is the sum of dihedral
angles of $n$-simplices containing the point.

(Notice that cone-manifolds whose singular locus has vertices are allowed.)

\begin{lemma}
\cite[pg.7]{BGP} If all cone angles of $n$-dimensional cone-manifold $M$ are
$\le 2\pi$, and $K \ge -1$, then $M$ is an Alexandroff space with curvature
bounded below by $-1$.
\end{lemma}

An $n$-dimensional cone-manifold clearly has Hausdorff dimension $n$.
Therefore we have the following corollary.

\begin{corollary}
\label{cone2} Let $M$ be an $3$-dimensional cone-manifold with all cone angles
$\le 2\pi$ and $K \ge -1$.  Let $M_0$ be a closed hyperbolic $3$-manifold.  If
$M$ and $M_0$ are homeomorphic then
$$\text{Vol} (M) \ge \text{Vol} (M_0).$$
\end{corollary}

All of these corollaries have only slightly weaker generalizations to
any dimension greater than $2$. We refer the reader to \cite{St}.  
There are other results whose methods follow along these lines; in
particular we would like to point out the work relating to Einstein
metrics, the Gromov norm, and simplicial volume for manifolds and covers
carried out by A.  Sambusetti (e.g. \cite{Sam1,Sam2}). However, these
recent results have already been well surveyed in \cite{BCG5}.

\section{Cautionary examples}
\label{section:counterexamples}

Our proof of the Degree Theorem (Theorem~\ref{theorem:degree}) can be
viewed as a step towards the Entropy Rigidity Conjecture (Conjecture
~\ref{conjecture:entrig}), or at least the inequality of that
conjecture; to do this one ``only'' needs to prove the inequality
(~\ref{equation:degree}) with the lowest possible $C$, namely
$C=\left(\frac{h(g)}{h(g_{0})}\right)^n$.

While the value of $C$ which comes out of the proof of the inequality
(see \cite{CF2}) can be explicitly computed, finding the best $C$
seems much harder.  In fact, it soon became clear to us that
the barycenter method applied without regard to the types of measures in
$\mathcal{M}(Y)$ is not sufficiently precise to obtain the rigidity aspect of
the theory, or even to prove that the locally symmetric metric minimizes
entropy (not necessarily uniquely).  In this section we describe some
explicit examples which demonstrate some of the problems.

The Jacobian estimates on the Besson-Courtois-Gallot map $F$ are obtained by
bounding the right hand side of (~\ref{eq:lastequation}), that is,
bounding
\begin{equation*}
\frac{\left(\det \ds{\int}_{K}
    O_{\theta}\begin{pmatrix} 1 & 0 \\
      0 & 0 \end{pmatrix}O_\theta^* \
    d\sigma_y(\theta)\right)^{\frac12}}{\det\ \ds{\int}_{K}
  O_{\theta}\begin{pmatrix} 0 & 0 \\ 0 & I
   \end{pmatrix}O_\theta^* \ d\sigma_y(\theta)}
\end{equation*}

independently of the measure.  In other words, no special property of
the measure $d\sigma_y$ is used.  This is also the case in the proof of
Entropy Rigidity in rank one (see \cite{BCG2}).

We recall that in a symmetric space the Furstenberg boundary is the
space of asymptotic (maximal) Weyl chambers which can be identified
with $G/P$ for a minimal parabolic $P$. This in turn is naturally
homeomorphic to $K/M$ where $K$ is the maximal compact subgroup of $G$
and $M$ is the centralizer of $\mathfrak{a}$ in $K$. Since $K$ is
isomorphic to $K_p$, isotropy subgroup of $p\in X$, and $M$ is
isomorphic to $M_p$, the centralizer of a Weyl chamber in
$\mathcal{F}_p$, we may interpret the Fursteberg boundary
geometrically. Namely it is naturally identified to any $G$-orbit of a
regular point in $\partial X$. We will denote by $\partial_FX$ the
particular orbit of the point at infinity corresponding to the
normalized barycenter $b^+$ of positive root vectors. By the above
$\partial_F X$ projects to the unit tangent sphere $S_pX$ as the orbit
of $b^+$ by $K_p$.

Hence if we write $J(\mu)$ for the square of the right hand side of
the expression \eqref{eq:Jacobian} with $\sigma_y^s$ replaced by $\mu$,
 then we would like for it to be
uniquely maximized when $\mu$ is the projection of Haar measure on
$K$. Here we only care about maximizing over the space of $M_p$
invariant $\mu$ supported on $K_p$.

\subsection{An example in $\hyp^2\times\hyp^2$}
In the case of $X=\hyp^2\times \hyp^2$ the Furstenberg boundary $\partial_FX$
is simply the torus $S^1\times S^1$. We will parameterize this by $(e^{i
t_1},e^{it_2})$. Let $M$ be a compact quotient of $X$
and $f:M\to M$ a continuous map. We first consider the
case $f=\op{Id}$ and the simple two-parameter family of probability measures
$\mu_{a,b}$ for $a,b>0$ on $\partial_FX$ given by
$$d\mu_{a,b}=\frac{2 + b + 2\,a\,{\cos (\frac{\Mvariable{t1}}{2})}^4 +
  b\,\cos (\Mvariable{t2})}{\left( 8 + 3\,a + 4\,b \right) \,{\pi
    }^2}d\mu_0$$
where $\mu_0$ the unit Haar measure on the torus. Note that $\mu_{0,0}$ is
Haar measure. For this family, by integrating over the torus we can compute
the expression for $J_0(\mu_{a,b}):=J(\mu_{a,b})/J(\op{Haar})$ to be
$$J_0(\mu_{a,b})=\frac{4\,{\left( 5\,a^2 + 24\,a\,\left( 2
        + b \right) + 16\,{\left( 2 + b \right) }^2 \right) }^2}
{{\left( 8 + 3\,a + 4\,b \right) }^2\,\left( 5\,a + 8\,\left( 2 + b
    \right) \right) \,\left( 7\,a + 8\,\left( 2 + b \right) \right)
  }.$$

\begin{figure}[ht]
\begin{center}
\psfrag{J}{\hspace{-.3in}$J(\mu_{a,b})$}
\psfrag{H}{}\psfrag{m}{}\psfrag{L}{}
\scalebox{1}{\includegraphics{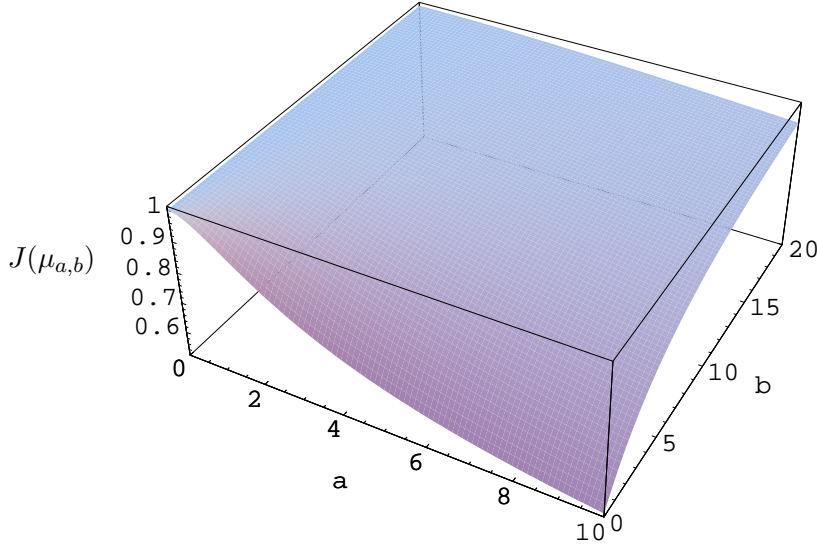}}
\end{center}
\caption{Graph of $J(\mu_{a,b})/J(\op{Haar})$ for several values of $a$
and $b$.}
\end{figure}

Note that the measures $\mu_{0,b}$ are distinct, but that
$J(\mu_{0,b})/J(\mu_{0,0})=1$ identically. In particular the Haar measure does
not uniquely minimize this quantity. This may not be that surprising in light
of the fact that $\hyp^2$ does not have a unique minimizer either nor is there
a entropy rigidity theorem for $X$, at least not for reducible lattices.

\subsection{An example in $\SL_3(\R)/\SO_3(\R)$}
Now we examine the case of $X=\op{SL}(3,\R)/\op{SO}(3)$. Suppose for a
single flat $\mathcal{F}_x$ and a sequence of $y_i\in \mathcal{F}_x$,
the measures $\sigma_{y_i}^s$ tend to the sum of Dirac measures
$\mu=\frac{1}{6}\sum_{i=1^6}\delta_{b_i^+(\infty)}$, where the $b_i^+$
are all the images of $b^+$ under the Weyl group in $\mathcal{F}_x$.
Hence there is one atom for each Weyl chamber at $x$ and they are
symmetrically placed making $x$ the barycenter of $\mu$. Then we claim
that if $N=M$ and $\tilde{f}$ induces the identity map on $\partial
\mathcal{F}_x$ then $\Jac F_s(y_i)\to\infty$.

The hypotheses that $N=M$ and $\tilde{f}$ is identity on $\partial
\mathcal{F}_x$, which contains the support of $\mu$, implies that
$\Jac F_s(y_i)$ is identical in the limit $i\to\infty$ to the estimate
on the right of \eqref{eq:Jacobian} when $\sigma_y^s$ is replaced by
$\mu$. Therefore it remains to show that this right hand side is
unbounded.

First note that the sum
$$\lim_{i\to \infty} \sum_{j=1}^6 dB_{(F_s(y_i),b_j^+(\infty))}^2$$ has
only a $3$-dimensional kernel, while
$$\lim_{i\to \infty}\sum_{j=1}^6 DdB_{(F_s(y_i),b_j^+(\infty))}$$ has
a $2$-dimensional kernel. Furthermore the numerator and denominator,
$$Q_1:=\int_{\D_F X} dB_{(F_s(y_i),\theta)}^2 d\sigma_{y_i }^s\ \
\mbox{and}\ \ Q_2:=\int_{\D_F X} DdB_{(F_s(y_i),\theta)}
d\sigma_{y_i}^s,$$
degenerate homogeneously. In particular, the
quantity $\det(Q_1)/\det(Q_2)^2$ is unbounded. This can be easily
verified explicitly in the case of a sum of eight Dirac measures for
which both integrals are nonsingular degenerating to the sum of the
six Dirac measures given above. The same phenomenon can easily be seen
to occur whenever we have an $\hyp^2$ factor.

In fact in \cite{CF2} we show that this cannot happen in
$\op{SL}(3,\R)$ by improving a result of Savage. However, as
demonstrated above the only reason is that the the measures
$\sigma_{y_i}$ never degenerate to any measure such as $\mu$.

\subsection{An example in $\SL_4(\R)/\SO_4(\R)$}
In the case of $\op{SL}(4,\R)/\op{SO}(4)$, Theorem
~\ref{theorem:degree} holds and the quantity \eqref{eq:Jacobian} is
bounded independently of $\sigma_y^s$. Nevertheless, here we present a
smooth (with respect to Haar) probability measure $\mu$ on
$\op{SO}(4)/M$ such that $J_0(\mu)$ is strictly larger than $1$.  In
fact, $\mu$ will be very close to the Haar measure $\mu_0$ only
differing by the adition of a few very small and sharp ``thorns'' and
removal of a few ``drillings.'' What we mean by this will become
evident from the final construction.

Our procedure is to consider second derivatives of the Jacobian
$$f(t):=J_0\left(\frac{\mu_0+t\mu}{\|\mu_0+t\mu\|}\right)$$ with respect
to a parameter $t$ where $\mu$ is an arbitrary signed measure. In fact,
if we can find such a $\mu$ for which the second derivative at zero is
positive then we will show how one can obtain such an example where
$\mu$ is a smooth (with respect to $\mu_0$) positive probability
measure.

We can verify directly that in general,
$J(\mu_0)=\left(\frac{s^2/n}{(\op{Tr}(DdB)/n)^2}\right)^n$ which
for $\op{SL}(4,\R)$ works out to be $\left(\frac{9
s^2}{40}\right)^9$. Hence, after clearing the denominator of the
measures, we can clearly rewrite
$$f(t)=\left(\frac{40}{9}\|\mu_0+t\mu\|\right)^9\frac{\det\left(\int_{\D_F X}
    d{B}_{(F_s(y),\theta)}^2
    d[\mu_0+t\mu](\theta)\right)}{\det\left(\int_{\D_F X}
    Dd{B}_{(F_s(y),\theta)}(\cdot,\cdot)
    d[\mu_0+t\mu](\theta)\right)^2}$$

Taking the derivative of the log we have, $\partial_t f(t)=f(t)
\partial_t \log f(t).$ Hence, the second derivative is,
$$\partial_t\partial_tf(t)=f(t)\left[
\left(\partial_t \log f(t)\right)^2+\partial_t\partial_t\log f(t)\right]$$

For an invertible matrix $M$ we have the operator identity $\log \det
M=\op{Tr} \log M$ and $\partial_t \op{Tr} \log M=\op{Tr}\left( (\partial_t
M)M^{-1}\right)$. Applying these identities we directly obtain,
\begin{gather*}
\begin{split}
\partial_t \log f(t)=\op{Tr} \left(\int_{\D_F X}
    d{B}_{(F_s(y),\theta)}^2
    d\mu(\theta)\right)\left(\int_{\D_F X}
    d{B}_{(F_s(y),\theta)}^2
    d[\mu_0+t\mu](\theta)\right)^{-1} +\\
    9 \frac{\|\mu\|}{\|\mu_0+t\mu\|}- 2\op{Tr} \left(\int_{\D_F X}
    Dd{B}_{(F_s(y),\theta)}
    d\mu(\theta)\right)\left(\int_{\D_F X}
    Dd{B}_{(F_s(y),\theta)}
    d[\mu_0+t\mu](\theta)\right)^{-1}
\end{split}
\end{gather*}

\smallskip
\noindent
Similarly,
\begin{gather*}
\begin{split}
\partial_t \partial_t \log f(t)=-\op{Tr}\left( \left(\int_{\D_F X}
    d{B}_{(F_s(y),\theta)}^2
    d\mu(\theta)\right)\left(\int_{\D_F X}
    d{B}_{(F_s(y),\theta)}^2
    d[\mu_0+t\mu](\theta)\right)^{-1}\right)^2 -\\
    9 \left(\frac{\|\mu\|}{\|\mu_0+t\mu\|}\right)^2+ 2\op{Tr}\left( \left(\int_{\D_F X}
    Dd{B}_{(F_s(y),\theta)}
    d\mu(\theta)\right)\left(\int_{\D_F X}
    Dd{B}_{(F_s(y),\theta)}
    d[\mu_0+t\mu](\theta)\right)^{-1} \right)^2
\end{split}
\end{gather*}
Since $f(0)=1$ and
$$\int_{\D_F X}
    d{B}_{(F_s(y),\theta)}^2
    d\mu_0(\theta)=\frac{1}{9}\Id \quad \text{and} \quad \int_{\D_F X}
    Dd{B}_{(F_s(y),\theta)}
    d\mu_0(\theta)= \frac{\sqrt{20}}{9}\Id$$
we see that $\partial_t\vert_{t=0}f(t)=0$ and
\begin{gather*}
\begin{split}
\partial_t\partial_t\vert_{t=0}f(t)=-9^2\op{Tr}\left(\int_{\D_F X}
    d{B}_{(F_s(y),\theta)}^2
    d\mu(\theta)\right)^2 -
    9 \left(\int_{\D_F X}d\mu\right)^2+\\
     \frac{9^2}{10}\op{Tr}\left( \int_{\D_F X}
    Dd{B}_{(F_s(y),\theta)}d\mu(\theta)\right)^2.
\end{split}
\end{gather*}

Now we consider the Cartan splitting $\op{SL}(4,\R)/K=KP$ where $P$ are the
positive definite symmetric matrices with determinant $1$ and we identify $K$
with the stabilizer of $F_s(y)$. $P$ is naturally identified with the $9$
dimensional subspace $\mathfrak{p}$ of traceless symmetric matrices in
$\mathfrak{sl}(4,\R)$. We also have the natural representation of $K$ in
$\op{SO}(\mathfrak{p})$ induced by the action of $K$  on $P$ by conjugation.

By lifting $\partial_FX$ to $K$ we may treat the integration over
$K$ and it is convenient to parameterize $K$ by $\R^6\equiv
\mathfrak{s}$ the subalgebra of skew-symmetric matrices in
$\mathfrak{sl}(4,\R)$. We can choose a domain of integration
$D\subset\R^6$, and integrate over this space. For $\theta\in D$ let
$\op{O}(\theta)\in\op{SO}(\mathfrak{p})$ be the corresponding element. If
$A$ is the constant diagonal matrix with diagonal $(1,0,0,0,0,0,0,0,0)$
and $B$ is the diagonal matrix with diagonal
$(0,0,0,\sqrt{\frac{2}{5}},\sqrt{\frac{2}{5}},\sqrt{\frac{2}{5}},
2\sqrt{\frac{2}{5}},2\sqrt{\frac{2}{5}},3\sqrt{\frac{2}{5}}),$ then
relative to the appropriate choice of coordinates on $T_{F_s(y)}X$ we may
write $d{B}_{(F_s(y),\theta)}^2=O(\theta) A O(\theta)^*$ and
$Dd{B}_{(F_s(y),\theta)}=O(\theta) B O(\theta)^*.$ Note that the metric we
are using is $\frac{1}{\sqrt{2}}$ times the one most commonly used by
representation theorists.

With these notations, the expression above for
$\partial_t\partial_t\vert_{t=0}f(t)$ becomes $\int_D\int_D
q(\sigma,\tau)d\mu(\sigma)d\mu(\tau)$ where
$$q(\sigma,\tau)=-9^2\op{Tr}\left(O(\sigma) A O(\sigma)^*O(\tau) A O(\tau)^*\right) -
    9 +\frac{9^2}{10}\op{Tr}\left( O(\sigma) B O(\sigma)^*O(\tau) B
    O(\tau)^*\right).$$

In particular, we may treat $\partial_t\partial_t\vert_{t=0}f(t)$
as a symmetric $2$-form $\Omega(\cdot,\cdot)$ on the space of
signed measures.

In fact, the function $q$ is equivariant so that $q(\sigma,\tau)$
can be written as $g(\sigma^{-1}\tau):=q(1,\sigma^{-1}\tau)$ and
one can easily check that $g(\sigma)=g(\sigma^{-1})$.

Three representative graphs of $g(\sigma)$ over two parameter
subspaces spanned by (respectively) $e_1$ and $e_2$, $e_3$ and
$e_4$, and $e_5$ and $e_6$, are given below. We see that the
function is mostly negative but there are smaller positive regions
as well.

\begin{figure}[ht]
\begin{center}
\scalebox{1}{\includegraphics{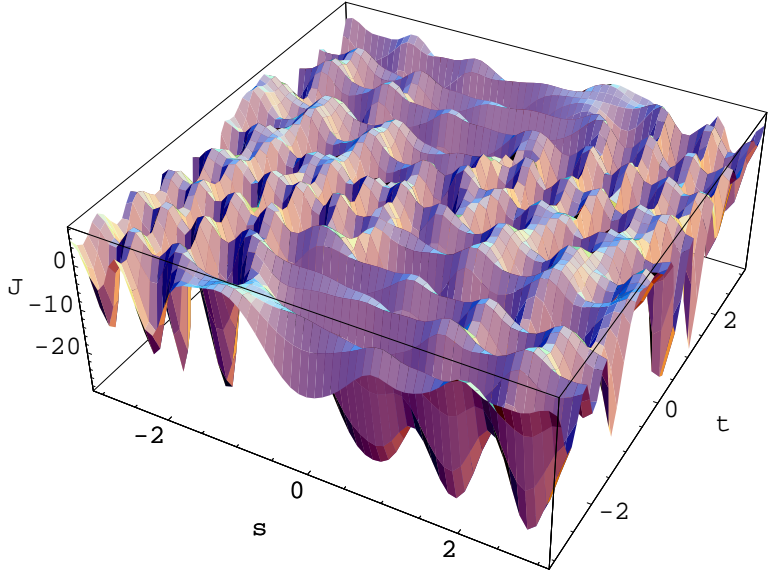}}
\end{center}
\caption{Graph of $g(se_1+te_2)$}
\end{figure}

\begin{figure}[ht]
\begin{center}
\scalebox{1}{\includegraphics{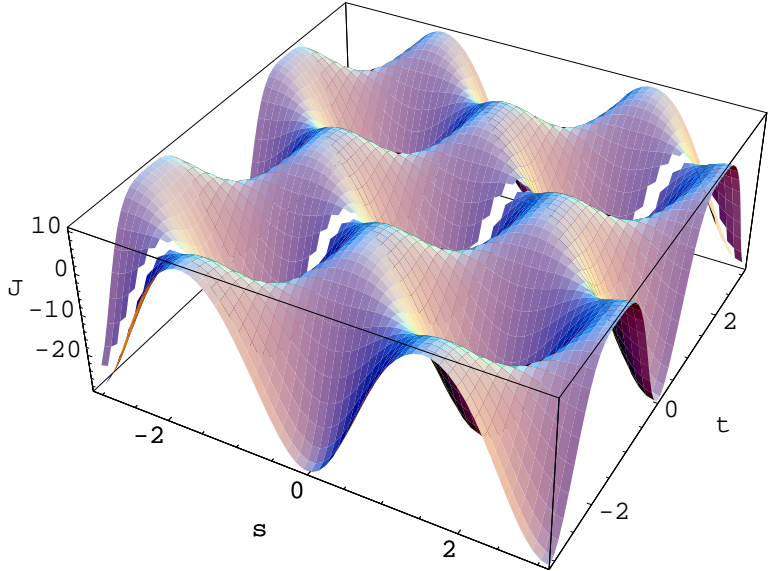}}
\end{center}
\caption{Graph of $g(se_3+te_4)$}
\end{figure}

\begin{figure}[ht]
\begin{center}
\scalebox{1}{\includegraphics{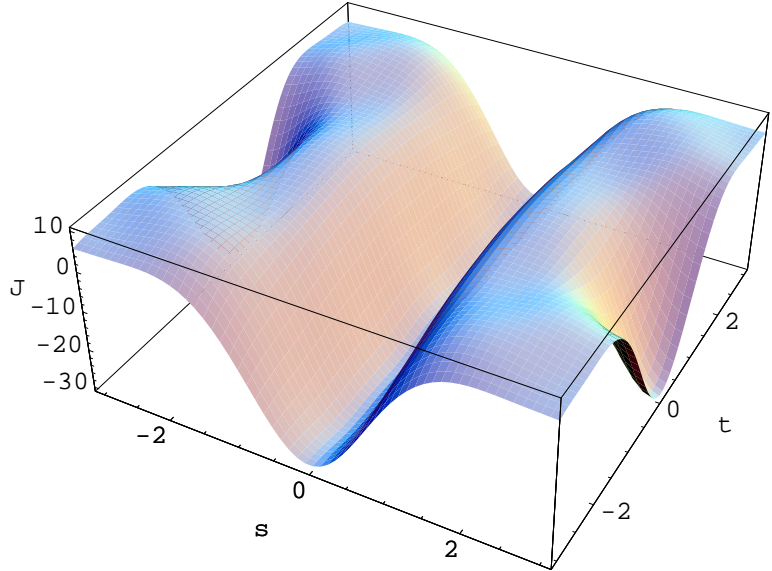}}
\end{center}
\caption{Graph of $g(s e_5+t e_6)$}
\end{figure}

Now we can take $\mu_1$ to be the sum of the atomic measures at
the $55$ points of $\partial_FX$ corresponding to the following
rational parameters in $D$. Namely the points are of the form
$K_i.b^+$ where $b^+$ is the barycenter of a set of positive roots
(any two choices will produce equivalent measures differing by a
pushforward by an isometry in the Weyl group) and $K_i$ is the
element of $K\eqsim SO(4)$ acting on $b^+$ by conjugation
corresponding to the element 
$$\exp
\begin{pmatrix}0 & s1 & s4 & s6 \\ -s1 & 0 & s2
& s5 \\ -s4 & -s2 & 0 & s3 \\ -s6 & -s5 & -s3 & 0
\end{pmatrix}\in SO(4)$$ 
For example, one such choice for $b^+$ has
coordinates $\{2/\sqrt{5},1/\sqrt{5},0\}$ in the canonical flat with the
most common choice of coordinates in $SL(4)$. The corresponding values of
the $s_1,\dots,s_6$ for $K_i$ may be read across the $i$-th row from the
following table: {\tiny
\begin{gather*}
\begin{matrix}
( -2.1763 & -1.1507 & -2.04369 & 0.44192 & 0.348006 & -0.0873793 )
\\ ( 1.10122 & -1.30303 & -1.07559 & 1.69668 & -0.537171 &
    -0.790495 ) \\ ( -2.66081 & -2.51729 & 0.874064 & 0.0170786 & 2.76331 & 0.918491 ) \\ ( 3.06666 & -1.20383 & 3.10876 & -1.8274 &
   1.33181 & 1.36971 ) \\ ( 2.14346 & 2.46489 & 0.233902 & -2.2138 & -1.34613 & -0.589322 ) \\ ( 2.27427 & 2.23082 & 2.87105 & 0.855587 &
    -0.33015 & -0.120273 ) \\ ( 2.39027 & 0.231281 & 1.93738 & 3.00424 & 2.76855 & 2.45438 ) \\ ( 2.01231 & 1.06647 & 2.80138 & 1.14019 &
    -2.46108 & 2.83835 ) \\ ( -2.48367 & 1.81689 & 0.446606 & 1.91056 & 2.00405 & -0.735379 ) \\ ( 1.31393 & 2.82133 & 2.2746 & 1.55063 &
    -1.49752 & -0.199991 ) \\ ( 3.02592 & -1.82225 & -0.293302 & -0.0626386 & -2.88422 & -1.13504 ) \\ ( 0.835977 & 2.01248 & -2.54401 &
   0.866367 & 0.15547 & 2.31572 ) \\ ( 3.08126 & 2.19107 & 2.85046 & -2.73643 & -2.06439 & -0.215145 ) \\ ( -1.60506 & -2.41617 & -1.19739 &
   1.37582 & 3.03405 & 0.925416 ) \\ ( -1.08172 & 0.0564745 & 0.185755 & -2.15354 & -1.3391 & -1.95008 ) \\ ( 2.49137 & -1.02443 &
    -1.93669 & 0.325143 & -0.805692 & -0.198555 ) \\ ( -1.87635 & 1.27567 & -0.514556 & -0.603716 & -2.95356 & -1.65078 ) \\ ( -2.05109 &
    -1.32914 & 1.38543 & 0.114991 & -1.94354 & 0.887035 ) \\ ( -0.674438 & -3.08308 & 1.0123 & -0.10102 & -2.47693 & 2.0086 ) \\ ( 1.66252 &
    -2.21819 & 2.60135 & -1.45814 & -0.67338 & 1.12196 ) \\ ( 1.33611 & 0.407789 & 2.98277 & -1.41592 & 1.14807 & -1.08302 ) \\ ( 1.89226 &
   3.05482 & 2.90424 & 1.94358 & 0.694209 & -0.973808 ) \\ ( 0.437082 & 1.88506 & 2.8235 & 2.26881 & -0.227581 & 3.01805 ) \\ ( -1.98061 &
   1.3454 & 0.31266 & 1.3346 & 1.83436 & -2.91815 ) \\ ( 2.11814 & -2.21478 & 1.99318 & 1.63935 & -2.17152 & 2.00983 ) \\ ( -3.04067 &
   1.72613 & -1.93417 & -3.07534 & -0.593287 & -0.441659 )
\end{matrix}
\end{gather*}
\begin{gather*}
\begin{matrix}
   ( 0.770342 & -1.81881 & -0.275197 & 0.431129 & -2.14367 & -1.69527 ) \\ (
    -1.43618 & 2.22732 & 0.685263 & 0.111728 & -0.128947 & 2.00388 ) \\ ( 1.70871 & -0.815082 & 1.01946 & -2.77706 & 0.738644 &
   0.316677 ) \\ ( 0.918539 & -1.36159 & -0.46878 & 0.250422 & -1.62977 & 2.22166 ) \\ ( 1.90247 & -1.07237 & 1.78702 & -1.35106 &
   0.904547 & -2.51869 ) \\ ( 0.0816078 & -0.436794 & -2.92231 & 0.511173 & -2.93104 & 0.700914 ) \\ ( -1.48943 & -1.81534 & -0.808906 &
   0.336383 & 0.913519 & 1.00958 ) \\ ( 1.41415 & -1.44362 & -1.75929 & -2.38244 & -0.0976769 & -0.52368 ) \\ ( -0.520171 & 1.83152 &
   1.25689 & -2.31421 & 1.71687 & 1.20862 ) \\ ( -1.96631 & 1.26418 & 1.49759 & -2.44414 & -2.17686 & -2.57833 ) \\ ( -0.154573 & 2.51279 &
   1.77364 & 0.226881 & 2.0735 & -1.63838 ) \\ ( -2.7821 & -1.4711 & 0.6912 & -2.39754 & 0.457165 & 2.19418 ) \\ ( -1.93022 & -1.08747 &
   2.34186 & 1.36679 & -0.505503 & 0.8455 ) \\ ( 1.16658 & -3.03898 & 1.1385 & 0.148052 & 0.201847 & 2.68094 ) \\ ( -1.84852 & 0.776858 &
   1.5698 & -0.687528 & -0.780427 & -0.72635 ) \\ ( 1.21031 & -2.35802 & 1.66997 & -1.4704 & -2.38844 & -1.41061 ) \\ ( 0.458594 & 2.75866 &
    -1.58872 & 0.364191 & -2.1775 & -1.22843 ) \\ ( 0.386299 & 0.261576 & -0.174403 & 1.76511 & -2.95714 & 0.722224 ) \\ ( -1.46748 &
    -2.15334 & -1.38535 & -1.73184 & 2.45454 & 1.7146 ) \\ ( 0.545927 & -2.51541 & -2.35701 & 0.0434106 & -0.207223 & 2.03679 ) \\ (
   0.325985 & 0.426342 & -1.7601 & -1.46899 & -0.638112 & -1.48682 ) \\ ( 0.995193 & 1.41102 & 2.67788 & -0.110336 & 0.810741 &
    -2.45279 ) \\ ( 1.00377 & -1.09859 & -0.9455 & 2.42064 & 1.69081 & 0.328404 ) \\ ( 1.65017 & 1.79446 & 0.906236 & -2.8566 & -1.2842 &
   2.89926 ) \\ ( -2.56134 & -0.141348 & -2.6657 & 1.22666 & 1.21836 & -1.79612 ) \\ ( -0.519297 & 2.95723 & 1.68207 & 1.45581 & 1.81155 &
   2.26842 ) \\ ( -2.46329 & -0.587199 & -0.384539 & 2.98937 & -1.01251 & 2.22599 ) \\ ( 1.10689 & -1.94668 & 1.22285 & 1.941 & -0.7505 &
    -1.70434 ) \\ ( 0.642599 & -1.05925 & -1.2264 & 0.210594 & 2.56583 & -2.40472
    ).
\end{matrix}
\end{gather*}}
Each of the $55$ atoms is weighted by the one of the three
following corresponding vector of weights: {\tiny
\begin{gather*}
\begin{split}
 \{& 0.0903801,0.0969632,0.145474,-0.0591319,-0.0214708,0.0919131,-0.000241965,0.0241533, \\
 & -0.0886868,0.104168,0.0392277,0.0888558,0.0227218,0.00893166,-0.212648,0.120529,\\
 &0.0324734,-0.0546992,0.111557,-0.159074,-0.130784,0.0239644,0.191044,   -0.212465,\\
 &-0.057288,-0.130236,-0.0515242,0.043965,0.0421474,0.0869512,0.0922342,0.186051,\\
 & 0.169673,0.365987,0.0227764,-0.0812135,-0.0426228,0.371221,-0.126781,-0.377116,\\
 &-0.172414,-0.00152986,-0.0305896,-0.00327764,-0.0494426,0.130664,-0.0487419,0.244147,\\
 &0.0402412,   -0.0658492,0.044671,-0.153327,-0.0747899,0.0137243,0.155563\}
    \\
   \{& 0.0237331,0.199178,0.123541,0.112583,-0.203652,-0.137257,0.11308,0.118312,\\
   &0.088158,-0.0285321,-0.00992443,0.106434,-0.0561753, 0.118035,-0.101555,0.10806,\\
   &-0.0693862,0.0678379,0.277438,-0.100565,-0.135501,0.105018,0.113,0.0720556,\\
   & -0.022502,0.116922,-0.0346328,0.0573696,-0.0847348,-0.257755,0.1629,-0.119937,\\
   &-0.0503877,-0.0524023,0.0808267,0.0327606,0.133811,-0.177716,0.0564386,0.0494853,\\
   & -0.075546,0.215228,-0.286515,-0.222407,-0.397166,0.081443,0.202924,-0.0419609,\\
   &-0.111385,-0.0710704,0.0650887,0.0620108,-0.150127,0.031396,0.0245637\} \\
    \{ &0.112757,0.172212,0.0460243,-0.0769192,0.151564,-0.0015713,0.280719,-0.194652,\\
    &0.0117092,-0.000721688,-0.177536,0.0200862,-0.0896234,-0.0198031,-0.0127717,-0.00823162,\\
    &0.0344373,-0.19736,0.137134,0.0857473,-0.00273806,-0.00872684,-0.105853,-0.157082,\\
    &0.0635683,0.106352,-0.13377,-0.304801,-0.0583263,0.220548,0.20856,0.0743173,\\
    &-0.184458,0.0193141,-0.000287821,-0.0463311,-0.223555,-0.284636,-0.0380433,0.00918501,\\
    &0.286329,0.131752,0.00101416,-0.0225734,-0.0120999,0.136267,-0.0511013,0.122597,\\
    &0.0157355, 0.183215,0.198623,-0.144791,0.0955358,0.0789536,0.157696\}
\end{split}
\end{gather*}}
For $\mu_1$ chosen to be any of these three weighted sums of
atomic measures on $\partial_F X$ it is not difficult (once the
representation for $K$ has been computed) to show that the
barycenter is $p$ in each case (the vector sum in $T_pX$ is $0$)
and $\Omega(\mu_1,\mu_1)$ is greater than $1.13346, 0.807823, $ or
$1.00141$ corresponding to the three sets of weights above.

The only problem is that $\mu_1$ is not a positive measure since each of
the above systems of weights has negative values. Nevertheless, we can fix
this by adding elements in the kernel of $\Omega$.  for any $\eps>0$ we
may take a smooth (w.r.t $\mu_0$) and symmetric (keeping barycenter $0$)
approximation $\mu_1^\eps$ to $\mu_1$ such that
$|\Omega(\mu_1,\mu_1)-\Omega(\mu_1^\eps,\mu_1^\eps)|<\eps$. Then if
$c(\eps)=\min\{\inf \frac{d\mu_1^\eps}{\mu_0},0\}$, the measure
$\mu_2=\frac{-c(\eps)\mu_0+\mu_1^\eps}{|c(\eps)|+\|\mu_1^\eps\|}$ is a
positive probability measure with barycenter $p$. Observe that $\int_D
O(\sigma)A O(\sigma)^* d\mu_0(\sigma)=\frac{1}{9}\Id$ and $\int_D
O(\sigma)B O(\sigma)^* d\mu_0(\sigma)=\frac{\sqrt{20}}{9}\Id$. It follows
from integrating inside the traces of the function $g(\sigma^{-1}\tau)$
that $\Omega(\mu_0,\mu)=0$ for any measure $\mu$. Therefore,
\begin{align}
\Omega(\mu_2,\mu_2)&=\frac{1}{(|c(\eps)|+\|\mu_1^\eps\|)^2}\Omega(\mu_1^\eps,\mu_1^\eps)\\
&\geq\frac{1}{(|c(\eps)|+\|\mu_1^\eps\|)^2} (0.807823-\eps)
\end{align}
which is positive for $\eps<0.8$.

The consequence of this is that $f(t)$ is larger than $1$ for some
sufficiently small $t$, and $\mu=\mu_2$. Hence there is no iso-derivative type
inequality for the Jacobian among (positive) probability measures which is
sharp at $\mu_0$.

In actuality the Jacobian at $\frac{\mu_0+t\mu_2}{\|\mu_0+t\mu_2\|}$ stays
less than $1.01$, and probably much less, because
$\frac{1}{(|c(\eps)|+\|\mu_1^\eps\|)^2}$ decays at least quadratically in
$\eps$ as $\eps\to 0$ for any choice of $\mu_1^\eps$.

\providecommand{\bysame}{\leavevmode\hbox to3em{\hrulefill}\thinspace}

\medskip

\noindent
Dept. of Mathematics, University of Chicago\\
5734 University Ave.\\
Chicago, Il 60637\\
E-mail: cconnell@math.uic.edu, farb@math.uchicago.edu

\end{document}